%% file: main.tex
\documentclass{IEEEtran} 
\usepackage{epsfig} 
\usepackage{mathptmx} 
\usepackage{times} 
\usepackage{graphics} 
\usepackage{epstopdf}
\usepackage{graphicx}      
\usepackage{hyperref}%
\usepackage{amsmath} 
\usepackage{amssymb}  
\newtheorem{remark}{Remark}
\newtheorem{definition}{Definition}

\usepackage{subfigure}
\newtheorem{theorem}{Theorem}

\usepackage{float}
\newtheorem{Assumption}{Assumption}
\newtheorem{proposition}{Proposition}
\graphicspath{ {images/} }
\usepackage{color}

\usepackage{enumerate}
\usepackage{caption}
\usepackage{tikz,pgfplots}
\usepackage{algpseudocode}
\def\RS{\textcolor{black}}

\usepackage[normalem]{ulem}

\begin{document}

\title{Safe Exploration and Escape Local Minima with Model Predictive Control under Partially Unknown Constraints}

\author{Raffaele Soloperto$^1$, Ali Mesbah$^{2}$, Frank Allg\"ower$^1$
\thanks{This work was supported by Deutsche Forschungsgemeinschaft (DFG, German Research Foundation) under Grant AL
316/12-2. Raffaele Soloperto thanks the International Max Planck Research School for Intelligent Systems (IMPRS-IS).}
\thanks{$^1$Institute for Systems Theory and Automatic Control, University of Stuttgart, Pfaffenwaldring 9, 70569 Stuttgart, Germany, e-mail: $\{$raffaele.soloperto, frank.allgower$\}$@ist.uni-stuttgart.de}
\thanks{$^2$Department of Chemical and Biomolecular Engineering, University of California, Berkeley, CA 94720, United States, e-mail:mesbah@berkeley.edu}}

\maketitle
\begin{abstract}
In this paper, we propose a novel model predictive control (MPC) framework for output tracking that deals with partially unknown constraints.
The MPC scheme optimizes over a learning and a backup trajectory. The learning trajectory aims to explore unknown and potentially unsafe areas, if and only if this might lead to a potential performance improvement. On the contrary, the backup trajectory lies in the known space, and is intended to ensure safety and convergence. 
The cost function for the learning trajectory is divided into a tracking and an offset cost, while the cost function for the backup trajectory is only marginally considered and only penalizes the offset cost.
We show that the proposed MPC scheme is not only able to safely explore the unknown constraints, but also escape from local minima that may arise from the presence of obstacles. Moreover, we provide formal guarantees for convergence and recursive feasibility of the MPC scheme, as well as closed-loop constraint satisfaction. Finally, the proposed MPC scheme is demonstrated in simulations using an example of autonomous vehicle driving in a partially unknown environment where unknown obstacles are present.
\end{abstract}



\input{Intro}
\input{problem_setup}

\input{preliminaries}
\input{proposed_approach}
\input{theoretical_analysis}

\input{example}

\section{Conclusion}
\label{sec:Conclusion}
In this paper, we presented a novel MPC framework for output tracking that deals with partially unknown constraints and ensures theoretical guarantees. 
We showed how the proposed MPC scheme optimizes over a learning and a backup trajectory, with the final goal to explore unknown areas while ensuring closed-loop safety.
The cost function for the learning trajectory is divided into a tracking cost, which minimizes the difference w.r.t. an online optimized setpoint, and an offset cost that minimizes the output of such a setpoint and the desired one. On the contrary, the cost function for the backup trajectory is only marginally considered and only minimizes the offset cost w.r.t. a potentially different setpoint.
We provided formal guarantees for convergence, recursive feasibility, and closed-loop constraint satisfaction. Finally, we showed the applicability of the proposed approach on a numerical example that considers an autonomous vehicle that drives in a partially unknown environment.

\bibliographystyle{ieeetran}  
\bibliography{Literature}

\begin{IEEEbiography}[{\includegraphics[width=1in,height=1.25in,clip,keepaspectratio]{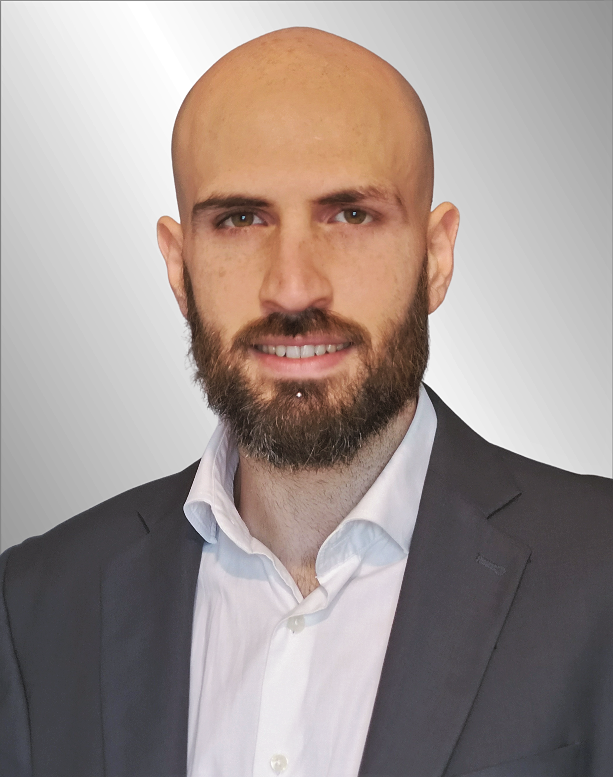}}]{Raffaele Soloperto}
 received both his bachelor's and master's degree in Automation Engineering from the University of Bologna, Italy, in June 2014 and October 2016, respectively.
Since July 2017, he is a Ph.D. student at the \emph{Institute for Systems Theory and Automatic Control} at the University of Stuttgart, Germany, and at the International Max-Planck Research School (IMPRS), under the supervision of Prof. Frank Allg\"ower. During his studies, he visited the TU-M\"unich, ETH-Z\"urich, EPF-Lousanne, and UC-Berkeley.
His research interests are in the area of robust, adaptive, and learning-based model predictive control.
\end{IEEEbiography}
\vspace{-3mm}
\begin{IEEEbiography}[{\includegraphics[width=1in,height=1.25in,clip,keepaspectratio]{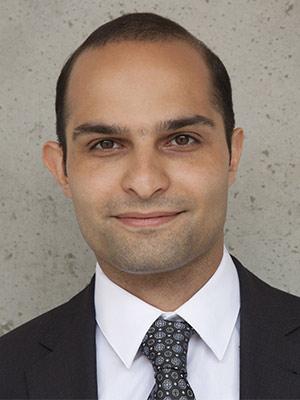}}]{Ali Mesbah} is Associate Professor of Chemical and Biomolecular Engineering at the University of California at Berkeley. Ali is a senior member of IEEE and AIChE. He is on the IEEE Control Systems Society conference editorial board as well as the editorial board of IEEE Transactions on Radiation and Plasma Medical Sciences and Optimal Control Applications and Methods. He is recipient of the Best Application Paper Award of the IFAC World Congress in 2020, the AIChE’s 35 Under 35 Award in 2017, the IEEE Control Systems Magazine Outstanding Paper Award in 2017, and the AIChE CAST W. David Smith, Jr. Graduation Publication Award in 2015. His research interests are in the area of learning and predictive control of uncertain systems.
\end{IEEEbiography}
\vspace{-3mm}
\begin{IEEEbiography}[{\includegraphics[width=1in,height=1.5in,clip,keepaspectratio]{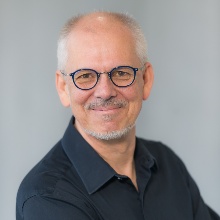}}]{Frank Allg\"ower}
 studied Engineering Cybernetics and Applied Mathematics in Stuttgart and at the University of California, Los Angeles (UCLA), respectively, and received his Ph.D. degree from the University of Stuttgart in Germany. Since 1999 he is the Director of the \emph{Institute for Systems Theory and Automatic Control} and professor at the University of Stuttgart.
His research interests include networked control, cooperative control, predictive control, and nonlinear control with application to a wide range of fields including systems biology.
For the years 2017-2020 Frank serves as President of the International Federation of Automatic Control (IFAC) and since 2012 as Vice President of the German Research Foundation DFG.
\end{IEEEbiography}

\end{document}

%% file: Intro.tex

\section{Introduction}
\subsubsection*{Motivation}
\begin{figure}[t]
\begin{center}
\includegraphics[width=0.47\textwidth]{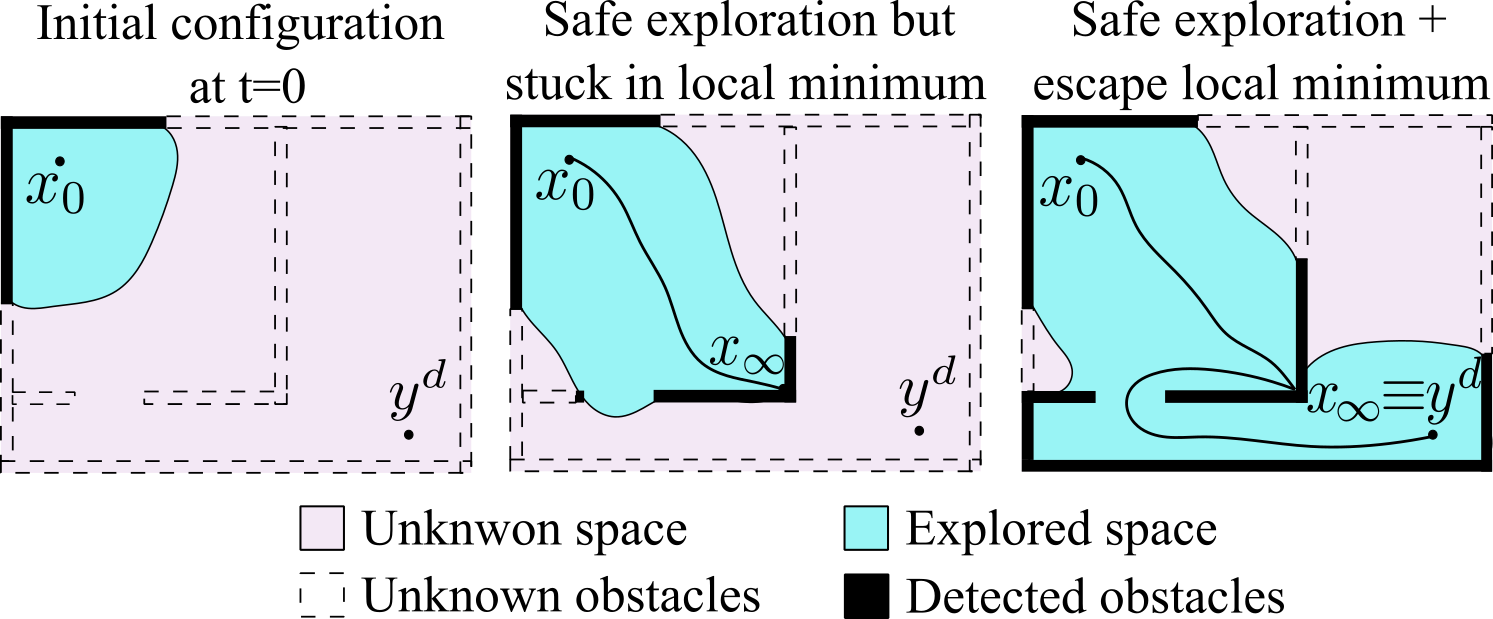}
\caption{Closed-loop representations of a system evolving in a partially unknown space under existing approaches (center) and proposed approach (right), with $x_{\infty} := \lim_{t\rightarrow \infty}x_t$. Goal: reach $y^d$ from $x_0$.}
\label{fig:intro}
\end{center}
\end{figure}

 Model predictive control (MPC)~\cite{rawlings2009model,kouvaritakis2016model}, \cite{grune2017nonlinear} is widely applied for optimization-based control of general constrained nonlinear systems.  
The performance of MPC largely depends on the accuracy of the prediction model, as well as on the accuracy of the constraint set.
Among others, the goal of Learning-based MPC approaches is to improve the closed-loop control performance by actively refining the model  (e.g., 
\cite{zanon2019safe, lorenzen2019robust,soloperto2020augmenting, kohler2019adaptive, solopertoguaranteed}) and the constraint set (e.g.,  \cite{tordesillas2021faster,saccani2021autonomous}) online. 
In several applications, such as autonomous driving, the system evolves in a partially unknown environment wherein the control objective is to get as close to a desired, but potentially non-reachable destination as possible. Due to the potentially changing environment wherein the system evolves, active exploration of unknown area may be necessary to improve the control performance. However, even though such an exploration process might be beneficial, it can also lead to violation of system constraints, loss of recursive feasibility of the MPC scheme, and instability of the system itself. Additionally, it is important to note that the presence of obstacles, combined with a (potentially non-convex) constraint set, might lead the system to converge into \emph{local minima}, as shown in Figure \ref{fig:intro}.


%
\subsubsection*{Related work}

Optimizing for performance while ensuring safety is the goal of many applications.
In theory, the ideal closed-loop performance under model or constraint inaccuracy can be realized using dual control~\cite{feldbaum1960dual}. However, except for a few specific problem setups (e.g., \cite{heirung2017dual,soloperto2019dual}), a computationally tractable reformulation of the dual control problem is, in general, not attainable. For this reason, several approaches modify the control problem by artificially introducing a heuristic cost function that explicitly excites the system, as well as novel learning formulations that have the goal to reduce uncertainties affecting the system, e.g., \cite{weiss2014robust,marafioti2014persistently}; see \cite{mesbah2018stochastic} for a recent survey on dual control. In \cite{tordesillas2019faster,tordesillas2021faster}, a control approach is proposed for known systems that evolve in a partially unknown environment and are equipped with a sensor that detects its surrounding. The approach relies on the notion of constructing two different trajectories, one for exploration and one for exploitation. The work additionally proposes a compelling computationally tractable solution for constructing convex safe sets, as demonstrated on several practical scenarios. On the other hand, the approach only focuses on navigation for UAVs, and does not provide theoretical guarantees, which is one of the objectives of this work. 
 In \cite{bonzanini2021perception,bonzanini2021stability}, the authors focus on an analogous problem setup, but additionally considered the case where the sensing quality is of statistical nature and can depend upon system operation. As such, a perception-aware chance-constrained MPC scheme is proposed that seeks to trade-off between the impact of control on sensing and sensing on control. 
The idea of constructing multiple trajectories is also exploited in \cite{lucia2014multi} with the goal of ensuring closed-loop constraint satisfaction in a robust MPC setting.
In \cite{saccani2021autonomous}, a multi-trajectory MPC approach is used for autonomous UAV navigation, wherein the problem setup only considers linear systems and does not provide convergence guarantees.

In order to increase the region of attraction of an MPC scheme, the approach in \cite{limon2008mpc} introduces the idea of online computation of a so-called \textit{artificial setpoint}. In particular, the cost function is divided into a tracking cost, which penalizes the distance between the system and the artificial setpoint, and an offset cost that steers the artificial setpoint towards the desired one.  The main advantages of the approach in \cite{limon2008mpc} are the increased region of attraction of the MPC scheme, the possibility to consider potentially unreachable setpoints, and ensuring recursive feasibility even when the desired state changes over time. This scheme has been extended for output tracking of nonlinear systems in \cite{limon2018nonlinear} and for linear systems with a mixture of deterministic and stochastic uncertainties in \cite{paulson2019mixed}. Recently, in \cite{solopertononlinear}, this scheme has been extended to the case of no terminal ingredients, which strongly simplifies its implementation. Finally, an approach for self-tuning terminal cost in economic MPC is analyzed in \cite{muller2013economic}. The paper proposes a novel idea to escape local minima by appropriately tuning the terminal cost online. However, the knowledge of the best achievable setpoint is required, and the overall problem setup differs from the one considered in this work since it does not consider learning or exploration.

\subsubsection*{Contributions}

For nonlinear systems with partially unknown constraints, this paper presents a novel MPC scheme that enables the controlled system to safely explore unknown areas within the state-space while ensuring convergence to the closest reachable setpoint. In particular, the MPC scheme considers two trajectories: a learning trajectory that is intended to improve the performance by exploring unknown areas of the state-space, and a backup trajectory to ensure safety as well as to find the best possible setpoint in the reachable space.
The novelty of the proposed MPC scheme is summarized as follows. 
\begin{itemize}
\item \textit{Optimal exploration:} The proposed MPC scheme mainly optimizes over the learning trajectory, and only \emph{marginally}  considers the backup trajectory by multiplying it with a small weight. This way the learning trajectory is left free to improve the performance by exploring unknown areas and, thus, does not involve a trade-off between exploration and exploitation.
\item \textit{Escaping local minima:} Thanks to the proposed cost function and constraints in the MPC scheme, we show that the system is able to escape, in closed-loop,  local minima that can possibly result from the presence of obstacles, or the given constraint set, as shown in Fig. \ref{fig:intro}.
\item \textit{Output tracking:} Motivated by practical applications and theoretical works, the proposed approach minimizes the distance between the output of the system and a desired output. The desired output can be non-reachable, and may change over time without influencing the theoretical properties of the proposed approach.
\item \textit{Theoretical guarantees:} In addition to recursive feasibility and constraint satisfaction, we propose an additional constraint in the MPC scheme to ensure convergence to the closest reachable setpoint. Using a simple motivating example, We demonstrate how the lack of such a constraint cannot ensure convergence.
\end{itemize}
As an additional contribution, we also extend the work presented in \cite{limon2008mpc} and in \cite{limon2018nonlinear}, to the case of non-connected steady-state manifolds. Such an extension is needed in order to have a more fair comparison between the theoretical results obtained in the proposed approach and the additional ones obtainable by extending \cite{limon2018nonlinear}.

\subsubsection*{Outline}
Section~\ref{sec:Problem_Formulation} presents the problem setup and shows the significance of exploring unknown areas within the state-space using a motivating example. In Section \ref{sec:preliminaries}, the preliminaries of the MPC scheme are presented, while discussing the shortcomings of an MPC scheme that employs only one trajectory for the problem setup at hand (i.e., partially unknown constraints).  
The proposed MPC scheme is presented in Section~\ref{sec:proposed_MPC_overall}, along with its theoretical analysis in \ref{sec:theoretical_analysis}. 
Section~\ref{sec:example} demonstrates the effectiveness of the MPC scheme using a numerical example and Section~\ref{sec:Conclusion} concludes the paper.
\color{black}

\subsubsection*{Notation}
 The quadratic norm with respect to a positive definite matrix $Q=Q^\top$ is denoted by $\|x\|_Q^2=x^\top Q x$.  
The minimal and maximal eigenvalue of $Q$ are denoted by $\lambda_{\min}(Q)$ and $\lambda_{\max}(Q)$, respectively. 
$\mathbb{R}_{\geq 0} = \{ r\in\mathbb{R}|r\geq  0\}$ denotes positive real numbers.
By $\mathcal{K}_{\infty}$, we denote the class of functions $\alpha:\mathbb{R}_{\geq 0}\rightarrow\mathbb{R}_{\geq 0}$, which are continuous, strictly increasing, unbounded and satisfy $\alpha(0)=0$. 
A variable $x$ predicted at time $t$ for $k$ steps ahead will be denoted by $x_{k|t}$, where $0\leq k\leq N$, with $N\in \mathbb{N}$ being the prediction horizon.  The notation $x_{\cdot|t}$ will be used when the entire predicted trajectory is considered. The variable $\mathbb{B}_a$ represents a hyper-ball in a given dimension with radius $a>0$, centered at the origin. 

%% file: problem_setup.tex

\section{Problem Setup}
\label{sec:Problem_Formulation}
In this paper, we consider the following nonlinear, time-invariant, discrete-time system
\begin{align}
\label{eq:perturbed_system}
{x}_{t+1} &= f(x_t, u_t), \nonumber \\
y_t &= h(x_t, u_t),
\end{align}
where $f: \mathbb{R}^n \times \mathbb{R}^m \rightarrow \mathbb{R}^n$ is the system dynamics, $h: \mathbb{R}^n \times \mathbb{R}^m \rightarrow \mathbb{R}^p$ is the ouptut function, $x\in \mathbb{R}^n$ is the state, $u\in \mathbb{R}^m$
is the control input, $y\in \mathbb{R}^p$ the output, and $t\in \mathbb{N}$ is the time.
At each time instant $t\geq 0$, the system is subject to the following state and input constraints
\begin{align}
\label{eq:constraints}
(x_t, u_t) &\in \mathbb{Z}\subseteq \mathbb{R}^{n+ m},
\end{align}
where $\mathbb{Z}$ is a compact and connected set.
Note that we do not require convexity of the set $\mathbb{Z}$ since this might directly include potential obstacles that the system must avoid.

\begin{figure*}[t!]
   \centering
\begin{tabular}{ccc}
\includegraphics[width=0.31\textwidth]{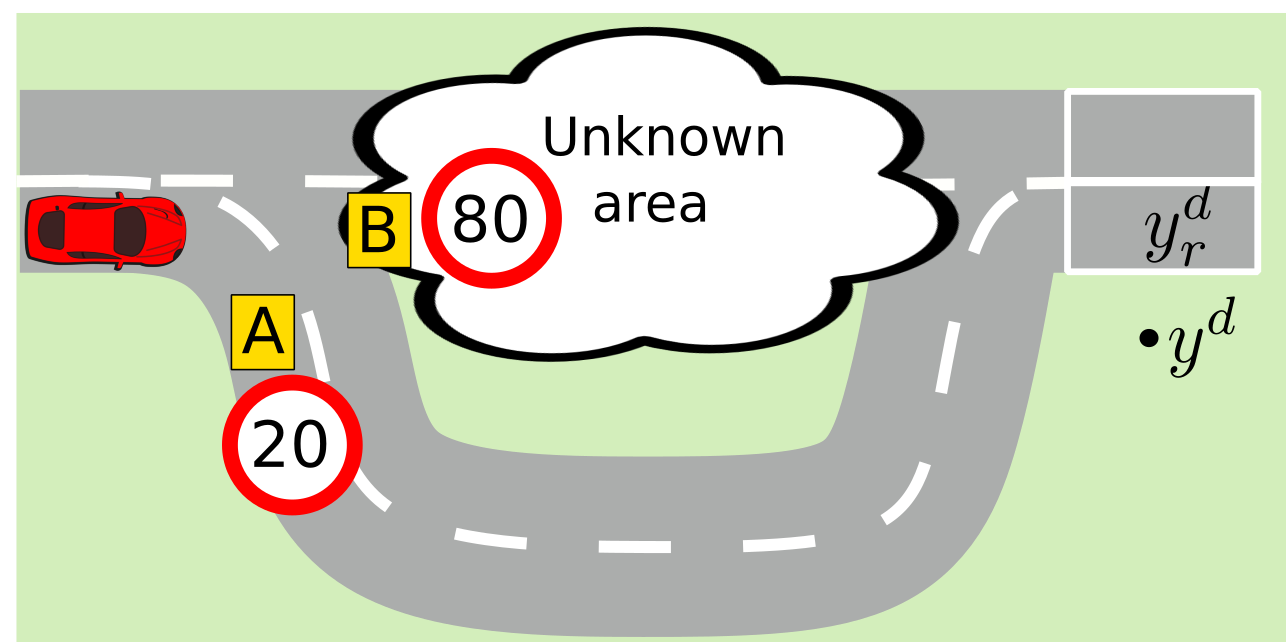}&
\includegraphics[width=0.31\textwidth]{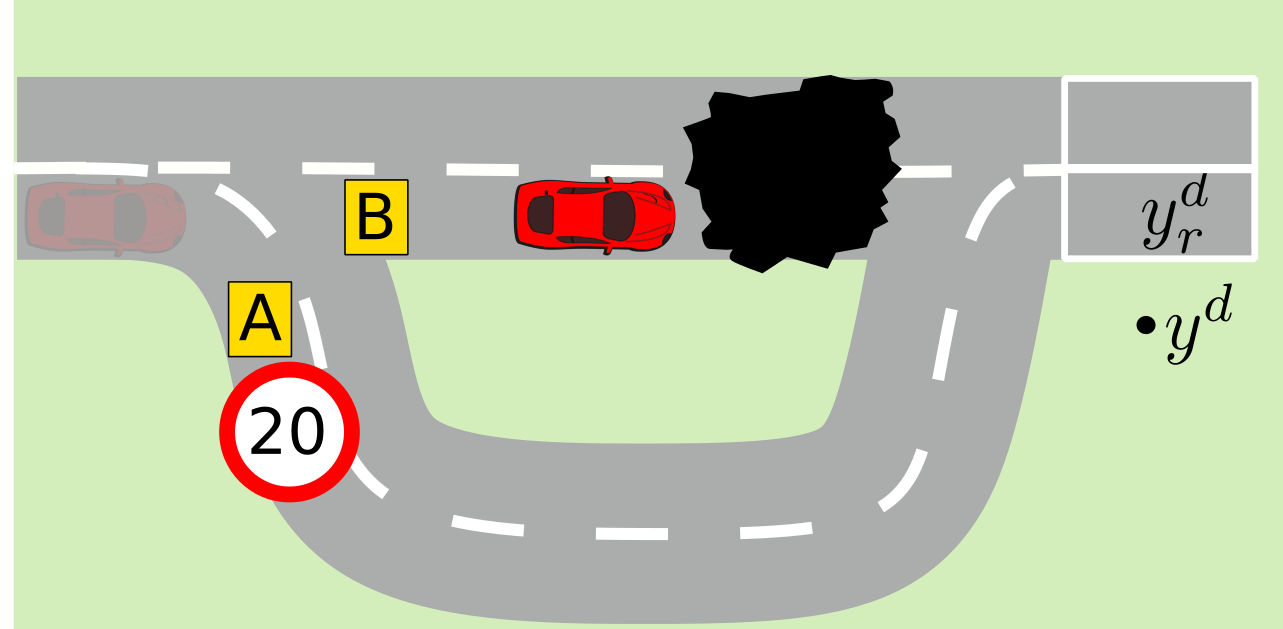} &
\includegraphics[width=0.31\textwidth]{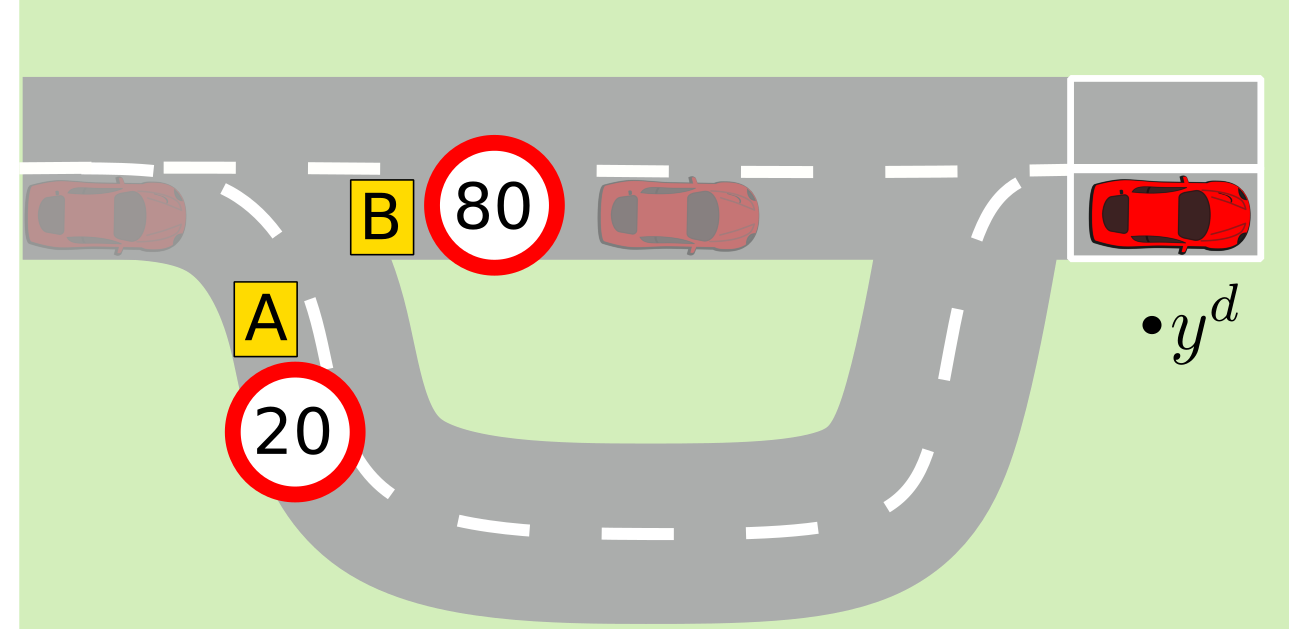}\\
Setup at $t=0$ & Scenario i & Scenario ii
\end{tabular}
   \caption{Trade-off between exploration and exploitation, and presence of a local minimum in Scenario i.}
\label{fig:motivation}
\end{figure*} 

We consider the case where the constraint set $\mathbb{Z}$ is only partially known, i.e., at each time instant $t\geq 0$, the set $\mathbb{R}^{n+ m}$ is subdivided into a known compact safe set $\mathbb{E}_t$ that satisfies
\begin{align*}
\mathbb{E}_t \subseteq \mathbb{Z},
\end{align*}
and a remaining unknown set $\mathbb{Z} \backslash  \mathbb{E}_t$.
\RS{\begin{remark}
Even though the presence of uncertainty in the input constraint set is rather limited in practical applications, we nevertheless consider it in order to include a wide variety of applications. This also implies that the more common case where only part of the state constraint is subject to uncertainty is also included in this problem setup.
\end{remark}}
\begin{Assumption}
\label{ass:neighborhood}
There exists a neighborhood $\mathbb{H}:\mathbb{R}^{n+m}\rightarrow 2^{\mathbb{R}^{n+m}}$ of the state and input pair $z_t=(x_t,u_t)$ that satisfies the following
\begin{align}
\label{eq:neighborhood}
[ z_t \oplus \mathbb{B}_h ] \cap \mathbb{Z}\subseteq \mathbb{H}(z_t) \subseteq \mathbb{E}_t, \quad \forall t\geq 0,
\end{align}
where $h>0$.
Moreover, at each time instant $t\geq 0$, an estimated constraint set $\mathbb{Z}_t$ is available, such that the following condition holds for all $t\geq 0$
\begin{align}
\label{eq:estimated_explored}
\mathbb{E}_t \subseteq \mathbb{Z}_t.
\end{align}
\end{Assumption}
As an example, in several autonomous system, the set $\mathbb{H}(z_t)$ might correspond to the set detected from a LIDAR sensor.
Note that the set $\mathbb{Z}_t=\mathbb{R}^{n+ m}$ trivially satisfies condition \eqref{eq:estimated_explored}. However, as it will become clear later, having an estimated constraint set $\mathbb{Z}_t$ as close as possible to the actual constraint set $\mathbb{Z}$ might improve the performance of the proposed approach.

\paragraph*{Steady-state manifold}
Given a set $\mathbb{A}\subseteq\mathbb{Z}$, the steady-state and -input manifold $\mathbb{S}(\mathbb{A})$ is defined as follows
\begin{align}
\label{eq:manifold}
&\mathbb{S}(\mathbb{A}):=\{(x,u)\in \mathbb{R}^{n+m}| (x,u) \oplus \mathbb{B}_{\lambda} \subseteq \mathbb{A}, ~ x =f(x,u)\},
\end{align}
where $\lambda\ll h$. Note that the set $\mathbb{B}_{\lambda}$ is only necessary in order to ensure that $\mathbb{S}(\mathbb{A})\subseteq \text{int}(\mathbb{A})$, which is required for the local controllability argument discussed in Assumption \ref{ass:bound_cost_function}, and later used in the theoretical analysis in Section \ref{sec:theoretical_analysis}.
\begin{remark}
Even if the set $\mathbb{A}$ is convex, the manifold $\mathbb{S}(\mathbb{A})$ can nevertheless be non-connected. In this paper, we explicitly deal with such a problem setup, and show the theoretical results that can be obtained with existing approaches and compare them with the proposed MPC framework.
\end{remark}

\paragraph*{Control goal}
\label{sec:control_goal}
The goal of the proposed approach is to steer the output $y_t$ of system  \eqref{eq:perturbed_system} as close as possible to a user-defined, potentially non-reachable,  output $y^d\in \mathbb{R}^p$, while satisfying the constraint set $\mathbb{Z}$ in \eqref{eq:constraints}. Simultaneously, exploration of the unknown and potentially unsafe space $\mathbb{Z}\backslash  \mathbb{E}_t$ must be incentivize only if this might lead to a performance improvement.

\subsection{Motivating example}
\label{sec:motivating_example}

Consider the case of an autonomous vehicle driving in a partially unknown environment, with the goal to get as close as possible to a desired and non-reachable output $y^d$ while satisfying a constraint set $\mathbb{Z}$, as shown in the left in Figure \ref{fig:motivation}. Note that since the output $y^d$ is not reachable, the only achievable goal is to reach the closest reachable output denoted with $y^d_r$. Such an output $y^d_r$ should, ideally, be directly computed by the MPC scheme, and therefore might not be known a priori.
Note that for the sake of simplicity, in the following example we focus on conveying the high-level motivations of the proposed approach, without employing rigorous mathematical tools.

Firstly, the car is subject to internal constraints that limit its velocity, acceleration, and steering angle, and are encoded in the set $\mathbb{Z}_{sc}$.
At time $t=0$, the safe set $\mathbb{E}_0$ is formed by the internal constraints $\mathbb{Z}_{sc}$, and the safe path ``A", i.e.,
\begin{align*}
\mathbb{E}_0 := \mathbb{Z}_{sc} \cap A.
\end{align*}
On the contrary, the path ``B" is unknown, and therefore cannot be considered as \textit{safe}. In order to enable the system to explore the path ``B", we define the set $\mathbb{Z}_0$ as follows
\begin{align*}
\mathbb{Z}_0 := \mathbb{Z}_{sc} \cap (A \cup B).
\end{align*}

%

In scenario i and scenario ii, we show two of the potential outcomes resulting from exploring the path ``B". In the following, we give a more detailed explanation about the trade-off between \textit{exploration} and \textit{exploitation}.

\textit{Exploitation:} As classical in the tracking MPC literature, \cite{rawlings2009model}, \cite{kouvaritakis2016model}, \cite{grune2017nonlinear}, \cite{limon2008mpc}, the MPC scheme exploits the safe constraint set $\mathbb{E}_t$ in order to find the optimal open-loop trajectory that minimizes a given cost function while ensuring constraint satisfaction. In the considered example, this means that the vehicle neglects the path ``B", and directly follows the path ``A", which is, however, potentially slower. 

\textit{Exploration:} 
Exploration can be achieved by actively incentivizing the system to explore the potentially unsafe set $\mathbb{Z}\backslash \mathbb{E}_t$. In the considered example, this means that the vehicle is directed towards the path ``B", even though this is not known at time $t=0$. 
As the system explores new safe areas, the safe set $\mathbb{E}_t$ expands, and the resulting performance might significantly improve. 
However, as shown in scenario i in Figure \ref{fig:motivation}, exploration of unknown areas might lead to constraint violation if not appropriately done, loss of recursive feasibility, as well as the possibility to get stuck in a new local minimum, as shown in Figure \ref{fig:intro}.
On the contrary, in scenario ii, path ``B" is clear and the vehicle is able to make use of the less stringent speed limit, and therefore arrives at the closest reachable output in a shorter amount of time.

In conclusion, the desired approach should not only incentivize exploration, but must simultaneously ensure constraint satisfaction, convergence, and should also escape from potential local minima.
Note that in Section \ref{sec:counter_example}, we consider the example above at a more mathematical level, where we show additional problems that existing approaches suffer from.


\color{black}

%% file: preliminaries.tex

\section{Preliminaries for MPC scheme}
\label{sec:preliminaries}
In this section, we discuss the preliminaries needed for the proposed MPC scheme, and we extend to the case of non-connected manifolds $\mathbb{S}(\mathbb{Z})$ the work done in \cite{limon2018nonlinear}, \cite{kohlernonlinear18} and \cite{solopertononlinear}. In particular, in Section \ref{sec:existing_MPC_scheme}, we firstly introduce an MPC scheme that employs the concept of artificial setpoints. Such a scheme will serve as a basis for the proposed MPC framework that will be later introduced and discussed in Section \ref{sec:proposed_MPC_overall}. In Section \ref{sec:tracking_cost}, we detail the tracking cost function as well as the terminal ingredients, while in Section \ref{sec:offset_cost} we discuss the so-called offset cost. In Section \ref{sec:shortcomings}, we present the theoretical analysis of such an MPC scheme extended to the case of non-connected manifolds, and we discuss its shortcomings.

\subsection{MPC scheme with artificial setpoints and non-connected manifold}
\label{sec:existing_MPC_scheme}
In the following,  we show an MPC scheme that is based on the concept of \textit{artificial setpoints} and that deals with a non-connected manifold $\mathbb{S}(\mathbb{Z})$. In a classical MPC scheme, see, e.g., \cite{rawlings2009model}, the optimization problem considers a user-defined tracking cost function which minimizes the distance between the predicted trajectory and a desired, reachable, setpoint $r^d=(x^d, u^d)\in \mathbb{S}(\mathbb{Z})$. On the contrary, the following MPC scheme simultaneously minimizes the tracking cost between the predicted trajectory and an online optimized artificial setpoint $r^s$, as well as an offset cost that penalizes the distance of the output of such an artificial setpoint $y^s=h(r^s)$, and the desired output $y^d$. This idea was initially introduced in \cite{limon2008mpc} for the case of linear systems, and is extended in \cite{limon2018nonlinear}, \cite{kohler2020nonlinear}, \cite{solopertononlinear} to the case of nonlinear systems, output tracking, and output tracking for MPC without terminal ingredients, respectively.
The main benefits of introducing an artificial setpoint are the followings:
\begin{enumerate}
\item The desired state and input $(x^d, u^d)$ that satisfy $y^d=h(x^d, u^d)$ do not need to be known, nor to exist. Indeed, feasibility of the MPC scheme is not influenced by the desired output $y^d$, which can therefore be even unreachable, i.e., $(x^d, u^d)\notin \mathbb{Z}$, or it might not correspond to a steady-output of \eqref{eq:perturbed_system}, i.e., $(x^d, u^d)\notin \mathbb{S}(\mathbb{Z})$.
\item The region of attraction of the MPC scheme increases compared to the one of a classical MPC scheme, i.e., where only the tracking cost function is considered. In particular, while in a classical MPC scheme the initial state $x_0$ needs to be in a \textit{neighborhood} of the desired state $x^d$, in the considered scheme the state $x_0$ only needs to be in a \textit{neighborhood} of the feasible manifold $\mathbb{S}(\mathbb{Z})$.
\item Initial feasibility of the MPC scheme ensures recursive feasibility of the scheme even when the desired output $y^d$ is subject to online changes. Such a property is not guaranteed in a classical MPC scheme and is crucial in several applications where stabilizing the origin is not the ultimate goal.
\end{enumerate}

At each time $t\geq0$, given the state $x_t$ and a compact set $\mathbb{E}_t \subseteq \mathbb{Z}$, the following optimization problem is solved:
\begin{subequations}
\label{eq:standard_MPC_scheme}
\begin{align}
\label{eq:MPC_cost_function}
\min_{x_{\cdot|t},u_{\cdot|t},r_{t}^s} ~ & \epsilon V_N(x_{\cdot|t},u_{\cdot|t}, r^{s}_t) + T(r^s_t)\\
s.t. \quad & x_{k+1|t} = f(x_{k|t},u_{k|t}),\\
\label{eq:MPC_existing_constraint}
&(x_{k|t}, u_{k|t})  \in \mathbb{E}_t,\\
\label{eq:terminal_region}
& x_{N|t} \in \mathbb{X}_f(r^s_t, \mathbb{E}_t),\\\label{eq:MPC_terminal_steady_state}
& r^s_t \in \mathbb{S}(\mathbb{E}_t),\\
& x_{0|t} = x_t,\\
& k = 0, \dots, N-1 \nonumber,
\end{align}
\end{subequations}
where $\epsilon>0$ is a user-defined parameter, $V_N$ is the tracking cost function that penalizes the distance between the open-loop trajectory $(x_{\cdot|t}, u_{\cdot|t})$ and an online optimized \textit{artificial setpoint} $r^s_t$, $\mathbb{X}_f$ is the terminal region, while $T$ is the offset cost function that penalizes the distance between the output of such an artificial setpoint $y^s_t = h(r^s_t)$ and the user-defined desired output $y^d$. For additional details on $V_N$ and $T$ compare \ref{sec:tracking_cost} and \ref{sec:offset_cost}, respectively.
\begin{remark}
Note that the constant $\epsilon$ is usually implicitly included in the function $V_N$, and therefore its ideal value is $\epsilon=1$. However, in this paper, we show the effects that $\epsilon$ plays on the MPC scheme \eqref{eq:standard_MPC_scheme}, and for this reason we explicitly include it.
\end{remark}
The optimal open loop trajectory resulting from \eqref{eq:standard_MPC_scheme} is denoted with $(x_{\cdot|t}^*, ~u_{\cdot|t}^*)$, and the optimal setpoint is referred as $r^{s,*}_t:=(x^{s,*}_t, u^{s,*}_t)$. Moreover, we define the function $V^*_N$
as follows
\begin{align*}
V^*_N(x_t, r^{s,*}_t, \mathbb{A}) = V_N(x_{\cdot|t}^*, u_{\cdot|t}^*, r^{s,*}_t).
\end{align*}
Finally, the resulting closed-loop is given as follows
\begin{align*}
x_{t+1} = f(x_t,u_t), \quad u_t:= u_{0|t}^*.
\end{align*}

\begin{remark}
The optimal trajectory $x_{\cdot|t}^*, u_{\cdot|t}^*, r^{s,*}_t$ resulting from \eqref{eq:standard_MPC_scheme} is analogous to the optimal trajectory resulting from a classical MPC scheme where only the tracking cost $V_N$ is minimized with respect to the optimal setpoint $r^{s,*}_t\in \mathbb{S}(\mathbb{Z})$ (and not to the desired output $y^d$), i.e., where the offset cost $T$ is not considered and the setpoint $r^{s,*}_t$ is given and not optimized. In Section \ref{sec:theoretical_analysis}, we will make use of this property to analyze the theoretical properties of the proposed framework.
\end{remark}

\begin{Assumption}
\label{ass:update_sets_standard}
The safe set $\mathbb{E}_t$ is updated so that the following holds
\begin{align}
\label{eq:set_update_standard}
(x^{*}_{k|t-1},u^{*}_{k|t-1}) & \in \mathbb{E}_t, \quad \forall k \in [0, \dots, N-1], \nonumber\\
\mathbb{X}_f(r^{s,*}_{t-1}, \mathbb{E}_{t-1}) & \subseteq \mathbb{E}_t, \quad \forall t\geq 1.
\end{align}
\end{Assumption}
\RS{
\begin{remark}
Note that if the computational power is especially limited, then the set $\mathbb{E}_t$ can be recomputed every several time instants since condition \eqref{eq:set_update_standard} trivially holds for $\mathbb{E}_t=\mathbb{E}_{t-1}$. Additionally, one can make use of the approach in \cite{tordesillas2021faster} where a computationally tractable way for online computation of convex sets in a non-convex environment is proposed.
\end{remark}}

\subsection{Tracking cost function $V_N$}
\label{sec:tracking_cost}
Given a user-defined prediction horizon $N \in \mathbb{N}$, and an admissible setpoint $r^s = (x^s, u^s) \in \mathbb{S}(\mathbb{Z})$, the tracking cost function $V_N$ is defined as follows
\begin{align}
\label{eq:cost_function}
V_N(x_{\cdot|t},  u_{\cdot|t}, r^s ) := \sum_{k=0}^{N-1} \ell(x_{k|t}, u_{k|t}, r^s) + V_f(x_{N|t}, r^s), 
\end{align}
where $(x_{\cdot|t}, u_{\cdot|t})\in \mathbb{Z}^{N+1}$ is a state and input sequence, $\ell:\mathbb{Z}\times \mathbb{S}(\mathbb{Z})\rightarrow \mathbb{R}_{\geq 0}$ is a user-defined stage cost, and $V_f:\mathbb{R}^n\times \mathbb{S}(\mathbb{Z})\rightarrow \mathbb{R}_{\geq 0}$ is the so-called terminal cost, detailed below.
The stage cost is designed so that the following assumptions hold.
\begin{Assumption}
\label{ass:stage_cost}
Given a compact set $\mathbb{A}\subseteq \mathbb{Z}$, there exist constants $a_1, a_2>0$ such that for all $(x,u)\in \mathbb{A}$ and $r^s=(x^s, u^s)\in \mathbb{S}(\mathbb{A})$, the following holds
\begin{align}
\label{eq:lower_uppper_stage_cost}
a_1 \|x-x^s\|^2 \leq \ell^*(x, r^s, \mathbb{A})\leq a_2\|x-x^s\|^2,
\end{align}
where $\ell^*$ is defined as follows
\begin{align}
\label{eq:minimum_stage_cost}
\ell^*(x, r^s, \mathbb{A}) := &\min_u \ell(x, u, r^s)\nonumber \\
&\text{s.t. } (x,u)\in \mathbb{A}.
\end{align}
Moreover, there exist constants $k_1^{\ell},k_2^{\ell}>0$ such that given two admissible setpoints $r_1^s, r_2^s\in \mathbb{S}(\mathbb{Z})$, and any pair $(x,u)\in \mathbb{Z}$, the following bound holds
\begin{align}
\label{eq:continuity_setpoint_stage}
\ell(x,u,r_1^s)\leq k_1^{\ell}\ell(x,u, r^s_2) + k_2^{\ell} \|r_1^s-r_2^s\|^2.
\end{align}
\end{Assumption}
In \cite{solopertononlinear}, it is shown how Assumption \ref{ass:stage_cost} is, for example, satisfied when the following standard quadratic stage cost is employed
\begin{align}
\label{eq:quadratic_stage_cost}
\ell(x,u, r^s) = \|x-x^s\|^2_Q + \|u-u^s\|^2_R,
\end{align}
where $Q,R\succ0$ are matrices of appropriate dimension.

\begin{Assumption}
\label{ass:bound_cost_function}
There exists a constant $\chi>0$ such that for all $\mathbb{A}\subseteq \mathbb{Z}$, and for all  $r^s=(x^s, u^s)\in \mathbb{S}(\mathbb{\mathbb{A}})$ then it holds that
\begin{align*}
\|x-x^s\|^2 \leq \chi \Rightarrow V^*_N(x,r^s, \mathbb{A}) \leq \gamma \ell^*(x, r^s, \mathbb{A}),
\end{align*}
where $\ell^*$ is defined in \eqref{eq:minimum_stage_cost}.
\end{Assumption}
Assumption \ref{ass:bound_cost_function} corresponds to a local exponential cost controllability
condition (cf. \cite{boccia2014stability} [Ass. 1]), and is trivially satisfied whenever (quadratically bounded) terminal costs can be computed offline, e.g., \cite{limon2018nonlinear}. Moreover, this assumption
holds if the linearization of the considered system around every admissible setpoint $r^s \in \mathbb{S}(\mathbb{A})$ is stabilizable, using the fact that $r^s\in \text{int}(\mathbb{A})$ (see \eqref{eq:manifold}) and $\ell$ quadratically bounded (Ass. \ref{ass:stage_cost}).
\subsubsection*{Terminal ingredients}
\label{sec:terminal_ingredients}
Given a set $\mathbb{A}\subseteq \mathbb{Z}$ and a setpoint $r^s\in \mathbb{S}(\mathbb{A})$, the terminal cost $V_f:\mathbb{R}^n\times\mathbb{S}(\mathbb{Z})\rightarrow \mathbb{R}_{\geq0}$, combined with a so-called terminal controller $\kappa_f:\mathbb{R}^n\times \mathbb{S}(\mathbb{Z})\rightarrow\mathbb{R}^m$, and a terminal region $\mathbb{X}_f(r, \mathbb{A})\subseteq \mathbb{A}$ are chosen so that the following condition holds for any $r^s = (x^s, u^s)\in \mathbb{S}(\mathbb{A})$, and any $x\in \mathbb{X}_f(r^s, \mathbb{A})$
\begin{align}
\label{eq:terminal_ingredients}
V_f(x^+, r^s) & \leq V_f(x, r^s)- \ell(x, \kappa_f(x, r^s), r^s),\nonumber \\
(x, \kappa_f(x, r)) &\in \mathbb{A}, \nonumber\\
x^+ &\in \mathbb{X}_f(r^s, \mathbb{A}),
\end{align}
with $x^+=f(x,u)$.

Note that especially for the case of nonlinear systems, designing a terminal region $\mathbb{X}_f$ that depends on a potentially time-varying set $\mathbb{A}$ can be especially challenging. For this reason, a more conservative, but easier solution is to employ a terminal equality constraint, i.e., 
\begin{align*}
\mathbb{X}_f(r^s, \mathbb{A}) := r^s,\quad \kappa_f(x, r^s) := u^s, \quad \forall r^s\in \mathbb{S}(\mathbb{A}).
\end{align*}
However,  depending on the considered problem setup, such a solution might reduce the region of attraction of the MPC scheme \eqref{eq:standard_MPC_scheme}. In order to overcome this problem, in \cite{solopertononlinear} it is shown how an MPC scheme that employs artificial setpoints can be also used without terminal ingredients. Even though in this paper we do not explicitly consider the case without terminal ingredients, we conjecture that the results from \cite{solopertononlinear} can be nevertheless applied in the proposed approach.

\subsection{Offset cost $T$}
\label{sec:offset_cost}
In the following, we discuss the technical conditions that the offset cost $T:\mathbb{S}(\mathbb{Z})\rightarrow\mathbb{R}_{\geq 0}$ must satisfy.
\begin{Assumption}
\label{ass:offset_cost}
There exist constants $k_0^T,k_1^T>0$, such that for any setpoint $r=(x,u)\in \text{int}(\mathbb{S}(\mathbb{Z}))$, there exists a setpoint $ \hat{r}= (\hat{x},\hat{u})\in \mathbb{S}(\mathbb{Z})$ and a function $\bar{\tau}:\mathbb{S}(\mathbb{Z}) \rightarrow \mathbb{R}_{>0}$ such that for any $\tau\in [0, \bar{\tau}(r)]$ the following holds 
\begin{subequations}
\begin{align}
\label{eq:bound_steady_state}
\|\hat{r}-r\| & \leq k_0^T \tau \|r-r^d_r(\mathbb{Z})\|,\\
\label{eq:bound_T_difference}
T(\hat{r}) - T(r)&\leq -k_1^T\tau\|r-r^d_r(\mathbb{Z})\|,
\end{align}
\end{subequations}
where $r^d_r(\mathbb{Z}):=(x^d_r(\mathbb{Z}),u^d_r(\mathbb{Z}))$ is the closest setpoint in the set $\mathbb{Z}$, defined as follows
\begin{align}
\label{eq:best_reachable_setpoint}
x^d_r(\mathbb{Z}),u^d_r(\mathbb{Z}) := \arg\min_{x,u} & ~ T(r) \\
\text{s.t. } & (x,u) \in \mathbb{S}(\mathbb{Z})\nonumber.
\end{align}
\end{Assumption}
Assumption \ref{ass:offset_cost} ensures that for any setpoint $r\in\text{int}(\mathbb{S}(\mathbb{Z}))$, we can find a new setpoint $\hat{r}\in\mathbb{S}(\mathbb{Z})$ which is closer to the desired output $y^d$, (cf.~\eqref{eq:bound_steady_state}) and has a smaller offset cost (cf.~\eqref{eq:bound_T_difference}).
It is crucial to note that Assumption \ref{ass:offset_cost} only considers setpoints that are strictly contained in the manifold $\mathbb{S}(\mathbb{Z})$, and therefore such a property does not need to hold for the setpoints that belong to the boundaries of $\mathbb{S}(\mathbb{Z})$. This implies that Assumption \ref{ass:offset_cost} can also be satisfied if the manifold $\mathbb{S}(\mathbb{Z})$ is non connected.

In \cite{solopertononlinear} it is shown how  the following offset cost
\begin{align*}
T(r) = \|y-y^d\|^2_P, \quad y=h(r),
\end{align*}
where $P$ is a positive definite matrix, satisfies Assumption \ref{ass:offset_cost} for the case of linear systems.

\RS{Now that the offset cost $T$ is formally defined, the desired control goal with respect to tracking described in Section \ref{sec:control_goal} can be mathematically reformulated as follows
\begin{align*}
\lim_{t \rightarrow \infty} y_t = h(r^d_r(\mathbb{Z})).
\end{align*}}

\subsection{Theoretical analysis of the MPC scheme \eqref{eq:standard_MPC_scheme}}

\begin{figure}[t]
\includegraphics[width=0.48\textwidth]{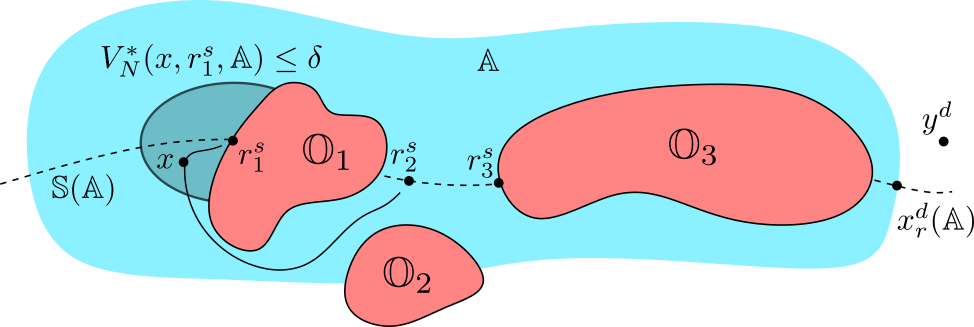}
\caption{Transitory setpoint ($r^s_1$) for the MPC scheme \eqref{eq:standard_MPC_scheme}, compare Definition \ref{def:transitory_steady_state_standard}. Obstacles are referred with $\mathbb{O}_i$, $i=1,2,3.$}\label{fig:transitory_setpoints}
\end{figure}

In the following, we show the theoretical analysis of the MPC scheme \eqref{eq:standard_MPC_scheme}, where we especially focus on the effect of $\epsilon$ on the closed-loop, and on the case of non-connected manifold $\mathbb{S}(\mathbb{Z})$. Even though this represents already an extension compared to existing works (\cite{limon2008mpc}, \cite{limon2018nonlinear}), it is only shown here with the purpose to have a fair comparison between the proposed approach in Section \ref{sec:proposed_MPC_overall}, and the one in \eqref{eq:standard_MPC_scheme}.

\begin{definition}
\label{def:transitory_steady_state_standard}
Given a set $\mathbb{A}\subseteq \mathbb{Z}$, we say that a setpoint $r^s_1\in \mathbb{S}(\mathbb{A})$ is a \textit{transitory setpoint} if there exists a scalar $\delta>0$ such that for all $(x,u)\in \mathbb{A}$ with $V^*_{N}(x, r^s_1, \mathbb{A})\leq \delta$ there exists an artificial setpoint $r^s_2\in \mathbb{S}(\mathbb{A})$ that satisfies the following
\begin{align}
\label{eq:local_setpoint_standard}
\epsilon V^*_{N}(x, r^s_2, \mathbb{A}) +T(r^s_2) <   \epsilon V^*_{N}(x, r^s_1, \mathbb{A}) + T(r^s_1).
\end{align}
We define with $\mathbb{G}(\epsilon, N, \mathbb{A})$ the set formed by all the transitory setpoints of the set $\mathbb{A}$, related to the scalar $\epsilon$ and to the MPC horizon $N$.
\end{definition}
The role that $\epsilon$ plays in Definition \ref{def:transitory_steady_state_standard} is discussed in Section \ref{sec:shortcomings}, where we discuss the shortcomings of the scheme \eqref{eq:standard_MPC_scheme}.
A graphical explanation of Definition \ref{def:transitory_steady_state_standard} is shown in Figure \ref{fig:transitory_setpoints}, where, in this case, the setpoint $r^s_1$ is a transitory setpoint, while the obstacle $\mathbb{O}_3$ is too large and therefore the setpoint $r^s_3$ is not a transitory setpoint according to Definition \ref{def:transitory_steady_state_standard}. It is important to note that Definition \ref{def:transitory_steady_state_standard} only states the sufficient conditions needed for overcoming an obstacle. This implies that it might be possible that the system will nevertheless overcome obstacle $\mathbb{O}_3$ and pass from $r^s_3$ to $r^d_r$ even if $r^s_3$ is not a transitory setpoint in the sense of Definition \ref{def:transitory_steady_state_standard}. 

Based on the definition of \textit{transitory setpoint}, in the following we propose a Theorem that is based on the result shown in \cite{solopertononlinear}, but that additionally includes the case of non-connected manifold as well as the set update $\mathbb{E}_t$.

\begin{theorem}
\label{thm:stability}
Let Assumptions \ref{ass:neighborhood}, \ref{ass:update_sets_standard},  \ref{ass:stage_cost}, \ref{ass:bound_cost_function}, and \ref{ass:offset_cost} hold, and suppose that at time $t=0$ Problem~\eqref{eq:standard_MPC_scheme} is feasible. Then, the MPC scheme is recursively feasible, and the system converges to a non-transitory setpoint $r^{s,*}_{\infty}$ that satisfies
\begin{align*}
r^{s,*}_{\infty} \in & \text{bound}(\mathbb{S}(\mathbb{Z}))\backslash \mathbb{G}(\epsilon,N, \mathbb{H} (r^{s,*}_{\infty})),\\
T(r^{s,*}_{\infty}) \leq & \lim_{t\rightarrow\infty} \inf_{i\in [0, t]} \epsilon V_N^*(x_i, r^{s,*}_i, \mathbb{H} (r^{s,*}_{\infty})) + T(r^{s,*}_i)\\
 \leq & \epsilon V_N^*(x_0, r^{s,*}_0, \mathbb{E}_0) + T(r^{s,*}_0).
\end{align*}
\end{theorem}
Theorem \ref{thm:stability} states that the system converges to a setpoint that belongs to the bounderies of the manifold $\mathbb{S}(\mathbb{Z})$, and that is not transitory for the set $\mathbb{G}(\epsilon,N, \mathbb{H} (r^{s,*}_{\infty}))$. This implies that if the $\mathbb{Z}$ is convex and $\mathbb{S}(\mathbb{Z})$ is connected, then the system converges to the closest reachable setpoint $r^d_r(\mathbb{Z})$, as defined in \eqref{eq:best_reachable_setpoint}, which is analogous to the result shown in \cite{kohlernonlinear18} and \cite{solopertononlinear}. On the contrary, if the manifold $\mathbb{S}(\mathbb{Z})$ is non-connected, then the system is nevertheless able to incrementally \textit{slide} on it, but it might converge to a non-transitory setpoint due to the potential presence of obstacles, or to the shape of the manifold itself.
The proof of Theorem \ref{thm:stability} is omitted since it can be obtained based on the proof of the proposed Theorem \ref{thm:convergence}, shown in Section \ref{sec:theoretical_analysis}.

\subsection{Shortcomings of the MPC scheme \eqref{eq:standard_MPC_scheme}}
\label{sec:shortcomings}
Even though the MPC scheme in \eqref{eq:standard_MPC_scheme} enjoys some interesting properties, such as convergence and constraint satisfaction, it is also affected by strong drawbacks that limit its applicability to a large class of problem setups. In the following, we list some of the major shortcomings that affect Problem \eqref{eq:standard_MPC_scheme}.

\paragraph*{Role of $\epsilon$}
Based on \eqref{eq:local_setpoint_standard}, it is easy to see that by decreasing the value of $\epsilon$, the MPC scheme \eqref{eq:standard_MPC_scheme} manages to overcome larger obstacles. In particular, given a state $x$ and two feasible setpoints $r_1^s, r_2^s\in \mathbb{S}(\mathbb{A})$, then the minimum value of $\epsilon$ required to go from $r_1^s$ to $r^s_2$ decreases as the cost $V^*_{N}(x, r^s_2, \mathbb{A})$ increases, according to the following condition
\begin{align*}
 V^*_{N}(x, r^s_2, \mathbb{A}) \stackrel{\eqref{eq:local_setpoint_standard}}{<}  \frac{T(r^s_1)-T(r^s_2)}{\epsilon} + V^*_{N}(x, r^s_1, \mathbb{A}),
\end{align*}
where $T(r^s_1)>T(r^s_2)$, and $V^*_{N}(x, r^s_2, \mathbb{A}) > V^*_{N}(x, r^s_1, \mathbb{A})$.
Note that for limit case of $\epsilon \ll 1$, we have that it is possible to overcome obstacles as long as $T(r^s_2)<T(r^s_1)$. However, it is nevertheless important to remark this is only possible if $V^*_{N}(x, r^s_2, \mathbb{A})$ exists, i.e., there must exist a feasible trajectory that starts in the $\delta$-neighborhood of the setpoint $r^s_1$, and stabilizes the setpoint $r_2^s$.

Even though choosing a very small value of $\epsilon$ is beneficial to overcome larger obstacles, it is easy to see that the cost function considered in the MPC scheme \eqref{eq:standard_MPC_scheme} behaves as follows
\begin{align}
\label{eq:max_obst_standard}
\epsilon\ll 1 ~ \Rightarrow ~\epsilon V_N(x_{\cdot},u_{\cdot}, r^{s}) + T(r^s) \approx T(r^s).
\end{align}
This means that, for $\epsilon\ll 1$, the MPC scheme \eqref{eq:standard_MPC_scheme} only tries to bring the output $y_t$ of the system as close as possible to the desired output $y^d$, while marginally (or not at all for $\epsilon \rightarrow 0$) considering the tracking cost $V_N$. Hence, since the transitory behavior of the system will not be optimized, the MPC scheme might make full use of the entire constraint set. Even though this might sound beneficial, note that this is undesired in several practical applications. For example, consider an autonomous vehicle controlled with the MPC scheme \eqref{eq:standard_MPC_scheme} with an $\epsilon\ll 1$. Then, the vehicle might travel at its maximum speed and acceleration, without considering fuel consumption, passenger comforts, nor the safety perceived by the surrounding systems, with the sole goal to get to the output $y^d$ as quick as possible. For this reason, in order to ensure that the system behaves according to a user-defined tracking cost function $V_N$, the value of $\epsilon$ cannot be chosen arbitrarily small, but must be appropriately tuned in order to find a trade-off between the maximum obstacle that can be overcome according to \eqref{eq:max_obst_standard}, and the desired behavior obtained with $\epsilon=1$, as discussed above.

\paragraph*{Lack of active exploration} The MPC scheme does not actively incentivize exploration of the unknown and potentially unsafe space $\mathbb{Z}\backslash \mathbb{E}_t$ since the open loop trajectory must entirely lie in the safe set $\mathbb{E}_t$, as enforced in condition \eqref{eq:MPC_existing_constraint}. This implies that, as also motivated in the motivating example in Section \ref{sec:motivating_example}, and shown in Figure \ref{fig:motivation}, there might exist safe trajectories with higher performance that will never be explored. The lack of active exploration is critical in several applications where the surrounding is not only unknown, but potentially also changing over time.

\color{black}

%% file: proposed_approach.tex

\section{Proposed MPC framework}
\label{sec:proposed_MPC_overall}

In order to overcome the shortcomings resulting from the MPC scheme \eqref{eq:standard_MPC_scheme} and discussed in Section \ref{sec:shortcomings}, in the following, we propose a novel MPC framework which differs from the one in \eqref{eq:standard_MPC_scheme} as follows
\begin{enumerate}
\item \textit{Overall idea:} Instead of optimizing only over one trajectory, the proposed MPC scheme optimizes over a so-called \textit{learning trajectory}, and a so-called \textit{backup trajectory}. At each time instant $t\geq 0$, both the trajectories share the current state $x_t$, and the first predicted input $u_{0|t}$ (therefore, the predicted state $x_{1|t}$ as well). The learning trajectory, denoted with the apex $L$, only satisfies the estimated constraint set $\mathbb{Z}_t$, and can therefore be unsafe since it might violate the actual constraint set $\mathbb{Z}$ (cf. Assumption \ref{ass:neighborhood}). On the contrary, the backup trajectory, denoted with the apex $B$, lies in the known constraint set $\mathbb{E}_t\subseteq \mathbb{Z}$.
In this way, we have that the learning trajectory is left free to explore the potentially unsafe set $\mathbb{Z}_t\backslash \mathbb{E}_t$, with the goal to find a trajectory with an improved open-loop performance, while the backup trajectory has the role to ensure safety as well as recursive feasibility of the MPC scheme.

\item \textit{Cost function: } In order to make sure that, in closed-loop, the system follows the learning trajectory rather than the backup one, the \textit{ideal} cost function should only minimize the learning trajectory, and neglect the backup one. In this way, the first predicted state $x_{1|t}^L=x_{1|t}^B$ is not computed based on a trade-off between the two trajectories, but it is only influenced by the learning trajectory, as long as there exists a feasible backup one. The \textit{ideal} cost function is formulated as follows
\begin{align}
\label{eq:ideal_cost_function}
V_N(x^L_{\cdot|t},u^L_{\cdot|t},r^{sL}_{t}) + T(r^{sL}_t).
\end{align}
In this case, the backup trajectory would not be optimized and simply results as an arbitrary feasible trajectory that lies in the set $\mathbb{E}_t$. However, as it will become clear in the theoretical analysis of the proposed scheme in Section \ref{sec:theoretical_analysis}, such an ideal cost function might not ensure convergence of the proposed MPC scheme. For this reason, we propose a minor modification of the cost function in \eqref{eq:ideal_cost_function} by additionally marginally minimizing the offset cost $T$ of the backup trajectory. Hence, the resulting cost function that will be employed in the proposed MPC scheme is as follows
\begin{align*}
V_N(x^L_{\cdot|t},u^L_{\cdot|t},r^{sL}_{t}) + T(r^{sL}_t)+ \epsilon T(r^{sB}_t),
\end{align*}
where, in this case, we have $\epsilon\ll 1$. 
\end{enumerate}

\subsection*{Proposed MPC scheme}

At each time $t\geq 0$, given the state $x_t$, the safe set $\mathbb{E}_t$ and the estimated constraint set $\mathbb{Z}_t$, the following optimization problem is solved:
\begin{subequations}
\label{eq:proposed_MPC}
\begin{align}
\min_{x^L_{\cdot|t},u^L_{\cdot|t},x^B_{\cdot|t},u^B_{\cdot|t},r^{sL}_{t},r^{sB}_{t}} & V_N(x^L_{\cdot|t},u^L_{\cdot|t},r^{sL}_{t}) +T(r^{sL}_t)+ \epsilon T(r^{sB}_t) \\
\label{eq:MPC_backup_dynamic}
\text{s.t.~} \quad & x^B_{k+1|t} = f(x^B_{k|t},u^B_{k|t})\\
\label{eq:MPC_backup_constraints}
\text{\small{Backup traj.:}}~ &(x^B_{k|t}, u^B_{k|t})  \in \mathbb{E}_t, \\
\label{eq:MPC_terminal_backup}
& x^{B}_{N|t} \in \mathbb{X}_f(x^{B,s}_{t},\mathbb{E}_t),~r^{B,s}_{t} \in \mathbb{S}(\mathbb{E}_t), \\[-6pt]
\rule{2cm}{0.5pt} & \nonumber \\ 
\label{eq:MPC_learning_dynamic}
& x^L_{k+1|t} = f(x^L_{k|t},u^L_{k|t},), \\
\label{eq:MPC_learning_constraints}
\text{\small{Learning traj.:}}~&(x^L_{k|t}, u^L_{k|t})  \in \mathbb{Z}_t,\\
\label{eq:MPC_terminal_learning}
& x^{L}_{N|t} \in \mathbb{X}_f(x^{L,s}_{t},\mathbb{Z}_t),~r^{L,s}_{t} \in \mathbb{S}(\mathbb{Z}_t),\\[-6pt]
\rule{2cm}{0.5pt} & \nonumber \\ 
\label{eq:MPC_initial_condition}
\text{\small{Initial cond.:}}~& x_{0|t}^L = x_{0|t}^B = x_t,  \\[-6pt]
\rule{2cm}{0.5pt} & \nonumber \\ 
\label{eq:MPC_first_input}
\text{\small{Safety:}}~&u^B_{0|t} = u^L_{0|t}, \\[-6pt]
\rule{2cm}{0.5pt} & \nonumber \\ 
\label{eq:sum_risk_negative}
\text{\small{Convergence:}}~& \epsilon V_N(x^B_{\cdot|t},u^B_{\cdot|t}, r^{sB}_t) + T(r_t^{sB}) \leq S_t + \hat{F}_t\\[-6pt]
\rule{2cm}{0.5pt} & \nonumber \\
& k = 0, \dots, N-1 \nonumber,
\end{align}
\end{subequations}
where condition \eqref{eq:sum_risk_negative} together with the variables $\hat{F}_t$ and $S_t$ are explained below. 

The optimal open-loop learning and backup trajectories obtained from \eqref{eq:proposed_MPC} are denoted with 
\begin{align*}
(x_{\cdot|t}^{L,*}, ~u_{\cdot|t}^{L,*}),  \quad r^{sL,*}_t:=(x^{sL,*}_t, u^{sL,*}_t),\\
(x_{\cdot|t}^{B,*}, ~u_{\cdot|t}^{B,*}),  \quad r^{sB,*}_t:=(x^{sB,*}_t, u^{sB,*}_t),
\end{align*}
while the resulting closed-loop is given as follows
\begin{align*}
x_{t+1} = f(x_t,u_t), \quad u_t:= u_{0|t}^{B,*}=u_{0|t}^{L,*}.
\end{align*}

Note that conditions \eqref{eq:MPC_backup_dynamic}-\eqref{eq:MPC_initial_condition} are analogous to the constraints used in the MPC scheme \eqref{eq:standard_MPC_scheme} but applied for the backup and learning trajectories.
Condition \eqref{eq:MPC_first_input} is included to ensure safety by forcing the system to remain in the safe set $\mathbb{E}_t\subseteq\mathbb{Z}$. Closed-loop exploration is then obtained thanks to the considered cost function, which will lead the system towards the unknown set $\mathbb{Z}_t\backslash\mathbb{E}_t$ if this improves the performance.
\paragraph*{Condition \eqref{eq:sum_risk_negative}}
In order to ensure convergence, we enforce, at each time step $t\geq 0$, a decrease of the backup cost function. For the sake of notational simplicity, we introduce the variable $F^{*}_t$ as follows
\begin{align}
\label{eq:F_optimal}
F^{*}_t:= \epsilon V_N(x_{\cdot|t}^{B,*}, ~u_{\cdot|t}^{B,*},~ r^{sB,*}_t) + T(r_t^{sB,*}).
\end{align}
Similarly, we define with $\hat{F}_{t+1}$ an upper-bound for the optimal value function $F^{*}_{t+1}$ computed at time $t$, and based on the optimal value function $F^{*}_{t}$, i.e.,
\begin{align}
\label{eq:F_backup}
\hat{F}_{t+1} := F^{*}_{t} 
- \epsilon\alpha \ell(x^{B,*}_{\cdot|t},u^{B,*}_{\cdot|t}, r^{sB,*}_{t}),
\end{align}
where $\alpha\in (0,1]$ is a user defined parameter.
Note that the value of $\alpha$ can be tuned to modify how often constraint \eqref{eq:sum_risk_negative} is active.
Starting from a user-defined $S_0\geq0$, the variable $S_t$ is updated as follow for all $t\geq 0$
\begin{align}
\label{eq:function_S_t}
S_{t+1} = & S_t -  F^{*}_t + \hat{F}_t.
\end{align}
Note that $S_t$ can be interpreted as a storage function.
\paragraph*{Initial feasibility of \eqref{eq:proposed_MPC}}
To ensure initial feasibility of the MPC scheme \eqref{eq:proposed_MPC} at time $t=0$, the value of $\hat{F}_0$ can be chosen arbitrarily large, and then $S_1=S_0$, or, alternatively, the constraint \eqref{eq:sum_risk_negative} can be implemented starting from $t=1$.
Moreover, we want to highlight that the proposed MPC scheme \eqref{eq:proposed_MPC} has the same region of attraction of the MPC scheme \eqref{eq:standard_MPC_scheme} since the open-loop learning trajectory can always be chosen as the backup one, ensuring satisfaction of the added conditions. Recursive feasibility is given in Theorem \ref{thm:convergence} and then shown in the theoretical analysis in Section \ref{sec:theoretical_analysis}.

\subsection{Motivating example for condition \eqref{eq:sum_risk_negative} }
\label{sec:counter_example}
\begin{figure}[t]
\begin{center}
\includegraphics[scale=0.4]{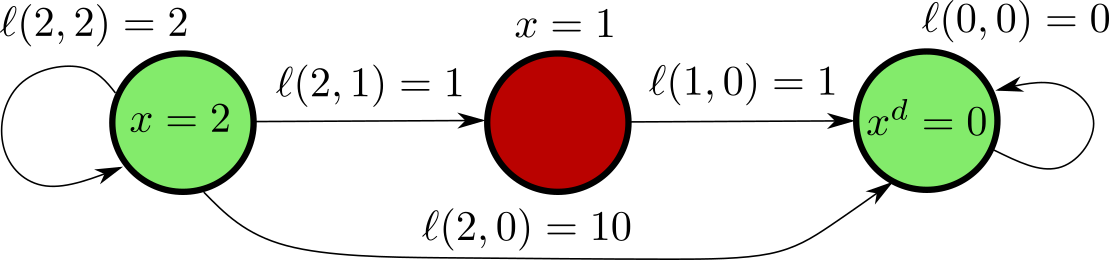}
\caption{The states $x_0=2$ and $x^d=0$ belong to the safe set, while the state $x=1$ does not. MPC horizon $N=3$.}
\label{fig:counter_example}
\end{center}
\end{figure}

In the following, we show through a motivating example, depicted in Figure \ref{fig:counter_example}, why condition \eqref{eq:sum_risk_negative} is needed to ensure convergence.
For the sake of simplicity, we consider a simplified version of the example discussed in Section \ref{sec:motivating_example}. By considering such a simple example, it is possible to clearly illustrate both the problems arising by the lack of condition \eqref{eq:sum_risk_negative}, as well as the effect that its presence has.
The considered constraint set $\mathbb{Z}$ is a non-convex subset of the set of natural numbers $\mathbb{N}^{n+m}$, which, therefore, contains only a finite number of points and hence the optimal trajectories can be easily computed.

Consider the following one-dimensional system 
\begin{align*}
x_{t+1}=u_t,
\end{align*}
with an MPC horizon $N=3$, initial state $x_0=2$, and desired state $x^d=0$. The sets $\mathbb{E}_t$, and $\mathbb{Z}_t$ are defined as follows
\begin{align*}
\mathbb{E}_t& :=\{(2,2), (2,0),(0,0)\},\\
\mathbb{Z}_t& :=\{(2,2),(2,1), (2,0),(1,0),(0,0)\},
\end{align*}
which means that the state $x=1$ is potentially unsafe, and therefore it belongs to the set $\mathbb{Z}_t$ and does not belong to the set $\mathbb{E}_t$. For the sake of simplicity, we enforce terminal equality constraint w.r.t the terminal region $\mathbb{X}_f := \{0\}$, and we do not consider the offset cost $T$. The stage cost $\ell(x,u)$ and the terminal cost $V_f$ are as follows
\begin{align*}
&\ell(2,2) = 2, \quad \ell(2,1) = 1, \quad \ell(2,0) = 10,\\
&\ell(1,0) = 1, \quad \ell(0,0) = 0, \quad V_f(0,0)=0, \end{align*}
and $\alpha=1$. 

In the following, we compare the case where condition \eqref{eq:sum_risk_negative} is not included, and the case where it is included in the MPC scheme \eqref{eq:proposed_MPC}.
Since condition \eqref{eq:sum_risk_negative} is not active at time $t=0$, we have that the optimal trajectories resulting from the MPC scheme \eqref{eq:proposed_MPC} with or without condition \eqref{eq:sum_risk_negative} are equivalent:
\begin{align}
\label{eq:optimal_example}
x_0 &= x_{0|0}^{L*} = x_{0|0}^{B*} =2,\nonumber \\
x_{1|0}^{B*} &=  x_{1|0}^{L*} = 2, \nonumber\\
x_{2|0}^{B*} &= 0,  \quad x_{2|0}^{L*} = 1,\nonumber \\
x_{3|0}^{B*} &= 0, \quad x_{3|0}^{L*} = 0,
\end{align}
which leads to $V_N(x_{\cdot|0}^{L*}, u_{\cdot|0}^{L*}) = 4$, and $V_N(x_{\cdot|0}^{B*}, u_{\cdot|0}^{B*}) = 12$. 

\begin{itemize}
\item Without ondition \eqref{eq:sum_risk_negative}:
At time $t=1$, we have that $x_1 = x_{1|0}^{B*} =x_0$. This implies that, since the trajectory \eqref{eq:optimal_example} remains feasible at time $t+1$ as well, the system evolves in closed loop as follows
\begin{align*}
x_t = 2, \quad \forall t\geq0.
\end{align*}
Therefore, even in such a simple scenario, and even in the presence of a safe trajectory, the system does not converge to the desired state $x^d=0$. Additionally, we want to remark that even by arbitrarily increasing the prediction horizon $N\geq 3$, then such a problem is not solved and the closed loop trajectory does not change.

\item With \eqref{eq:sum_risk_negative}:
Based on the optimal trajectory shown in \eqref{eq:optimal_example}, we have that, according to \eqref{eq:F_backup}, $\hat{F}_1$ is defined as follows
\begin{align*}
\hat{F}_1 = V_N(x_{\cdot|0}^{B*}, u_{\cdot|0}^{B*}) - \ell(x_{0|0}^{B*}, u_{0|0}^{B*}) = 12-1 = 11.
\end{align*}
Starting from a user-defined $S_0\geq0$, the value of $S_1$ is computed as follows
\begin{align*}
S_1 \stackrel{\eqref{eq:function_S_t}}{=}  S_0 - 12+11 = S_0-1.
\end{align*}
By iterating the reasoning above, we have that as long as the system remains at $x=2$, the value of $S_t$ is updated as $S_t = S_0 - t\cdot 1, \quad t\geq0$.
Therefore, this implies that, regardless of the initial value of $S_0$, there will be at time instant $t^*>0$ where $S_t + \hat{F}_t=11$, and hence the optimal solution shown in \eqref{eq:optimal_example} is not feasible anymore. This means that at time $t^*$ the only feasible trajectory is the following safe trajectory
\begin{align}
\label{eq:optimal_example_good}
x_{t^*} &= x_{0|t^*}^{L*} = x_{0|t^*}^{B*} =2,\nonumber \\
x_{1|t^*}^{B*} &=  x_{1|t^*}^{L*} = 0, \nonumber\\
x_{2|t^*}^{B*} &= 0,  \quad x_{2|t^*}^{L*} = 0,\nonumber \\
x_{3|t^*}^{B*} &= 0, \quad x_{3|t^*}^{L*} = 0,
\end{align}
which also implies $x_t=0$ for all $t\geq t^*+1$.
\end{itemize}
In conclusion, with such an example, we showed that even for simple problem setups, and for arbitrarily large prediction horizons $N\geq 3$, an MPC scheme that employs two trajectories without additionally including constraint \eqref{eq:sum_risk_negative} does not ensure convergence to the desired setpoint.

\color{black}


\begin{Assumption}
\label{ass:update_sets}
Similar to Assumption \ref{ass:update_sets_standard}, the safe set $\mathbb{E}_t$ is updated at each time instant $t\geq 0$ so that the following holds
\begin{align}
\label{eq:set_update}
(x^{B,*}_{\cdot|t-1},u^{B,*}_{\cdot|t-1}) & \in \mathbb{E}_t, \quad \forall k \in [0, \dots, N-1], \nonumber\\
\mathbb{X}_f(x^{sB,*}_{t-1}, \mathbb{E}_{t-1}) & \subseteq \mathbb{E}_t.
\end{align}
\end{Assumption}
Condition \eqref{eq:set_update} is trivially satisfied if the set $\mathbb{E}_t$ is chosen as $\mathbb{E}_t=\mathbb{E}_{t-1}$, which implies that the set $\mathbb{E}_t$ can be recomputed every several time instants, if the computational power is especially limited.
Note that the estimated constraint set $\mathbb{Z}_t$ is updated so that condition \eqref{eq:estimated_explored} (i.e., $\mathbb{E}_t\subseteq\mathbb{Z}_t$) holds also with the new safe set $\mathbb{E}_t$. 

%
%

\begin{figure}[t]
\includegraphics[width=0.48\textwidth]{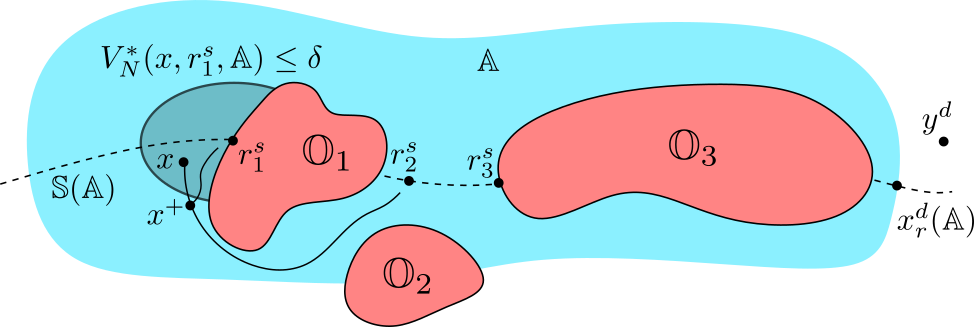}
\caption{Transitory setpoint ($r^s_1$) for the MPC scheme \eqref{eq:proposed_MPC}, compare Definition \ref{def:transitory_steady_state}. Obstacles are referred with $\mathbb{O}_i$, $i=1,2,3$.}\label{fig:transitory_setpoints_2}
\end{figure}
\begin{definition}
\label{def:transitory_steady_state}
Consider a set $\mathbb{A}\subseteq \mathbb{Z}$, we say that an artificial setpoint $r^s_1\in \mathbb{S}(\mathbb{A})$ is a \textit{transitory setpoint} for the MPC scheme \eqref{eq:proposed_MPC}, if there exists a scalar $\delta>0$ such that for all $(x,u)\in \mathbb{A}$ with $V_N^*(x, r^s_1,\mathbb{A}) \leq \delta$ there exists an artificial setpoint $r^s_2\in \mathbb{S}(\mathbb{A})$ with $T(r^s_2)<T(r^s_1)$ that satisfies the following
\begin{align}
\label{eq:transitory_setpoints_2}
& \epsilon\left[\ell(x,u, r^s_2) + V^*_{N-1}(x^+, r^s_2, \mathbb{A})\right] +T(r^s_2) \nonumber \\
 <  & \epsilon V_N^*(x, r^s_1,\mathbb{A}) + T(r^s_1),
\end{align}
where $x^+=f(x,u)$.
We define with $\mathbb{M}(\epsilon, N,\mathbb{A})$ the set constructed by all the transitory setpoints related to the set $\mathbb{A}$, the scalar $\epsilon>0$, and the MPC horizon $N$.
\end{definition}
In Figure \ref{fig:transitory_setpoints_2} we give a graphical representation of Definition \ref{def:transitory_steady_state}.

\begin{proposition}
\label{prop:sliding_steady_state}
Consider a set $\mathbb{A}\subseteq \mathbb{Z}$, and let Assumption \ref{ass:offset_cost} holds. If $r\in \text{int}(\mathbb{S}(\mathbb{A}))$, then $r$ is a transitory setpoint for the MPC scheme \eqref{eq:proposed_MPC}, i.e., 
$r\in \mathbb{M}(\epsilon, N, \mathbb{A})$, for any $\epsilon>0$.
\end{proposition}
The proposition above directly implies that a non-transitory setpoint $r$ must belong to the boundaries of the manifold $\mathbb{S}(\mathbb{A})$.
The proof of Proposition \ref{prop:sliding_steady_state} is shown in the Section \ref{sec:proof_prop_sliding_setpoint}.

\paragraph*{Role of $\epsilon$ in the MPC scheme \eqref{eq:proposed_MPC}}
Differently from the MPC scheme \eqref{eq:standard_MPC_scheme}, in the proposed MPC framework, the value of $\epsilon$ can be chosen arbitrarily small without yielding undesired effects. In particular, choosing a small $\epsilon$ leads to the following effects:
\begin{itemize}
\item The cost function optimized in the MPC scheme \eqref{eq:proposed_MPC} approaches the ideal cost function discussed in \eqref{eq:ideal_cost_function}, i.e.,
\begin{align*}
&\epsilon \RS{\ll} 1  \Rightarrow V_N(x^L_{\cdot|t},u^L_{\cdot|t},r^{sL}_{t}) +T(r^{sL}_t)+ \epsilon T(r^{sB}_t) \\ \RS{\approx} & V_N(x^L_{\cdot|t},u^L_{\cdot|t},r^{sL}_{t}) +T(r^{sL}_t).
\end{align*}
This implies that the closed-loop does not result from a trade-off between the learning and backup trajectory, but it focuses on the learning trajectory, leading to a safe exploration of more performing trajectories.
\item Based on Definition \ref{def:transitory_steady_state}, as the value of $\epsilon$ decreases, then it is possible to overcome larger obstacles, (i.e., escape \textit{deeper} local minima). In particular, given a cost function $V^*_{N-1}(x^+, r^s_2, \mathbb{A})$ needed to overcome an obstacle, then the following condition must be satisfied
\begin{align*}
& \ell(x,u, r^s_2) + V^*_{N-1}(x^+, r^s_2, \mathbb{A}) \\ < & \frac{T(r^s_1)-T(r^s_2)}{\epsilon} + V_N^*(x, r^s_1,\mathbb{A}).
\end{align*}
The condition above implies that if $T(r^s_1)>T(r^s_2)$, then it is always possible to choose a small enough $\epsilon$ to overcome an obstacle (assuming that both $V^*_{N-1}(x^+, r^s_2, \mathbb{A})$ and $V_N^*(x, r^s_1,\mathbb{A})$ exist). We want to remark again that choosing a small $\epsilon>0$ it is beneficial for the proposed MPC scheme \eqref{eq:proposed_MPC}, while for the MPC scheme \eqref{eq:standard_MPC_scheme} the value of $\epsilon$ must be chosen based on a trade-off between the size of the obstacle that the system can overcome and the desired performance.
\end{itemize}

%

\begin{theorem}
\label{thm:convergence}
Let Assumptions  \ref{ass:neighborhood}, \ref{ass:stage_cost}, \ref{ass:bound_cost_function}, \ref{ass:offset_cost}, and \ref{ass:update_sets} hold, and suppose that at time $t = 0$ Problem \eqref{eq:proposed_MPC} is feasible. Then, the MPC scheme is recursively feasible for all $t\geq 0$, and the following holds
\begin{align*}
\lim_{t\rightarrow \infty} z_t = r^{s,*}_{\infty},
\end{align*}
where $r^{s,*}_{\infty}$ satisfies the following
\begin{align}
\label{eq:manifold_theorem}
&r^{s,*}_{\infty}\in \text{bound}(\mathbb{S}(\mathbb{Z})) \backslash \mathbb{M}(\epsilon,N,\mathbb{H}(r^{s,*}_{\infty})), \quad \text{and}\\
\label{eq:decrease_theorem}
T(r^{s,*}_{\infty}) \leq & \lim_{t\rightarrow\infty} \inf_{i\in [0, t]} \epsilon  V_N(x_{\cdot|i}^{B,*}, ~u_{\cdot|i}^{B,*},~ r^{sB,*}_i) + T(r^{s,*}_i) + S_i\nonumber \\
 \leq & \epsilon  V_N(x_{\cdot|0}^{B,*}, ~u_{\cdot|0}^{B,*},~ r^{sB,*}_0)  + T(r^{s,*}_0) + S_0.
\end{align}
\end{theorem}
Note that the results shown in Theorem \ref{thm:stability} and in Theorem \ref{thm:convergence} are qualitatively analogous. However, the formal benefit resulting form Theorem \ref{thm:convergence} is the possibility to arbitrarily decrease the value of $\epsilon$ in order to increase the size of the set $\mathbb{M}(\epsilon,N,\mathbb{H}(r^{s,*}_{\infty}))$. This implies that the system is able to overcome larger obstacles, without deteriorating the performance, which is not possible in the case of the MPC scheme \eqref{eq:standard_MPC_scheme}. Additionally, we want to remark that the proposed MPC scheme \eqref{eq:proposed_MPC} enables the system to explore unknown areas, which represents a main contribution related to the scheme \eqref{eq:standard_MPC_scheme}.
The proof of Theorem \ref{thm:convergence} is shown in Section \ref{sec:theoretical_analysis}.


%% file: theoretical_analysis.tex

\section{Theoretical Analysis}
\label{sec:theoretical_analysis}

In this section, we show the proof of Proposition \ref{prop:sliding_steady_state} in  Section \ref{sec:proof_prop_sliding_setpoint}, 
and the proof of Theorem \ref{thm:convergence} in Section \ref{sec:recursive_feasibility} and \ref{sec:convergence}, where we first discuss recursive feasibility and then convergence, respectively.

\subsection{Proof of Proposition \ref{prop:sliding_steady_state}}
\label{sec:proof_prop_sliding_setpoint}
The proof of Proposition \ref{prop:sliding_steady_state} is done by contradiction.

Firstly, considering the artificial setpoint $r^s=(x^s, u^s)\in \text{int}(\mathbb{S}(\mathbb{A}))$ and based on Assumption \ref{ass:offset_cost}, we know that there exists a different artificial setpoint $\hat{r}^s=(\hat{x}^s,\hat{u}^s)\in \mathbb{S}(\mathbb{A})$ that satisfies the following
\begin{subequations}
\begin{align}
\label{eq:move_steady_state}
\|\hat{r}^s - r^s\| &\leq k_0^T \tau \|r^s-r^d_r(\mathbb{A})\|,  \\
\label{eq:evo_offset}
T(\hat{r}) - T(r) & \leq -k_1^T\tau\|r-r^d_r(\mathbb{Z})\|,
\end{align}
\end{subequations}
with a later specified value of $\tau\in (0,\bar{\tau}(r^s))$.

In the following, we show that for a small enough $\tau$, there exists a feasible trajectory that stabilizes the artificial setpoint $\hat{r}^s$ and that has a smaller cost compared to the trajectory that stabilizes the artificial setpoint $r^s$. Therefore, this results in a contradiction, and implies that the setpoint $r^s$ is a transitory setpoint in the sense of Definition \ref{def:transitory_steady_state}.

\paragraph*{Existence of a feasible solution for $\hat{r}^s$} 
In order to ensure that the MPC problem is feasible, we focus on the following neighborhood of the artificial setpoint $r^s=(x^s,u^s)$ 
\begin{align*}
\|x-x^s\|^2\leq \chi,
\end{align*}
which allows us to make use of Assumption \ref{ass:bound_cost_function}, which ensures the existence of a solution that stabilizies $r^s$.

Since we assume that $r^s \neq r^d_r(\mathbb{A})$, then it is always possible to  upper bound the optimal cost function as follows
\begin{align}
\label{eq:bound_stage_a_bar}
 V_N^*(x, r^s,\mathbb{A}) \leq  \bar{a}\|r^s-r^d_r(\mathbb{A})\|^2,
\end{align}
with a later specified value of $\bar{a}>0$.
Based on the properties of the stage cost, we can construct the following upper bound
\begin{align}
\label{eq:bound_norm_new_state}
& 
\|x^+-\hat{x}^s\|^2 = \|x^+-x^s+x^s-\hat{x}^s\|^2 \nonumber\\
= &\|x^+-x^s\|^2+\|x^s-\hat{x}^s\|^2 + 2\|x^+-x^s\|\|x^s-\hat{x}^s\| \nonumber\\
\stackrel{\eqref{eq:cost_function},\eqref{eq:lower_uppper_stage_cost}, \eqref{eq:move_steady_state},\eqref{eq:bound_stage_a_bar}}{\leq} & \frac{\bar{a}}{a_1}\|r^s-r^d_r(\mathbb{A})\|^2   + k_0^{T,2} \tau^2 \|r^s-r^d_r(\mathbb{A})\|^2 \nonumber \\
&+ k_0^T \tau \sqrt{\frac{\bar{a}}{a_1}}  \|r^s-r^d_r(\mathbb{A})\|^2.
\end{align}
Hence, one can choose $\bar{a}$ and $\tau$ small enough to ensure the following
\begin{align}
\label{eq:bound_state_new}
\|x^+-\hat{x}^s\|^2 \leq \chi, 
\end{align}
which allows us to invoke again Assumption \ref{ass:bound_cost_function} to ensure the existence of a feasible trajectory that starts in $x^+$ and stabilizes the setpoint $\hat{r}^s$ within $N-1$ steps. It is crucial to note that, in this proof, we do not change the first pair $(x,u)$, with $x^+=f(x,u)$, so that condition \eqref{eq:MPC_first_input} holds.

\paragraph*{Cost comparison between $r^s$ and $\hat{r}^s$}
In the following, we show that by stabilizing the setpoint $\hat{r}^s$ (without changing the first state and input pair $(x,u)$), the overall cost function is smaller compared to the one resulting from stabilizing the setpoint $r^s$, i.e., that the following holds
\begin{align}
\label{eq:desired_result_cost}
& \epsilon\left[\ell(x,u, \hat{r}^s) + V^*_{N-1}(x^+, \hat{r}^s, \mathbb{A})\right] +T(\hat{r}^s)\nonumber \\
< & \epsilon V_N^*(x, r^s,\mathbb{A}) + T(r^s).
\end{align}
Firstly, we analyze each \textit{unknown} term of \eqref{eq:desired_result_cost} individually, and then we combine the obtained bounds to show satisfaction of \eqref{eq:desired_result_cost}.
\begin{itemize}
\item Stage cost $\ell(x,u, \hat{r}^s)$:
\begin{align}
\label{eq:bound_stage_cost_new}
&\ell(x, u, \hat{r}^s)\stackrel{\eqref{eq:continuity_setpoint_stage}}{\leq} k_1^{\ell}\ell(x,u, r^s) + k_2^{\ell} \|r^s-\hat{r}^s\|^2 \nonumber \\
\stackrel{\eqref{eq:move_steady_state}, \eqref{eq:bound_stage_a_bar}}{\leq}& k_1^{\ell}\bar{a}\|r^s-r^d_r(\mathbb{A})\|^2 + k_2^{\ell} k_0^{T,2} \tau^2 \|r^s-r^d_r(\mathbb{A})\|^2,
\end{align}
where we used the fact that $\ell(x,u, r^s)\stackrel{\eqref{eq:cost_function}}{\leq} V^*_{N}(x,u, r^s)$.
\item Optimal cost function $V^*_{N-1}(x^+, \hat{r}^s, \mathbb{A})$:
\begin{align}
\label{eq:bound_short_cost}
&V^*_{N-1}(x^+, \hat{r}^s, \mathbb{A})\stackrel{\text{Ass.} \ref{ass:bound_cost_function}}{ \leq} \gamma\ell^*(x^+, \hat{r}^s, \mathbb{A}) \nonumber \\ \stackrel{\eqref{eq:lower_uppper_stage_cost}}{\leq} &\gamma a_2 \|x^+ - \hat{x}^s\|^2 \nonumber \\
\stackrel{\eqref{eq:bound_state_new}}{\leq} &  \gamma a_2 \left( \frac{\bar{a}}{a_1}+
k_0^{T,2} \tau^2 + k_0^T \tau \sqrt{\frac{\bar{a}}{a_1}}\right) \|r^s-r^d_r(\mathbb{A})\|^2.
\end{align}

%

\item Offset costs $T(\hat{r}^s)- T(r^s)$: Shown in \eqref{eq:evo_offset}.
\end{itemize}

Finally, we can choose $\bar{a}$ and $\tau$ small enough, and combine \eqref{eq:evo_offset}, \eqref{eq:bound_stage_cost_new}, and \eqref{eq:bound_short_cost}, to obtain the following bound

\begin{align}
\label{eq:final_prop}
&\epsilon\left[\ell(x_t,u_t, \hat{r}^s) + V^*_{N-1}(x^+, \hat{r}^s, \mathbb{A})\right]  \underbrace{-\epsilon V^*_N(x_t, r^s, \mathbb{A})}_{ \leq  0 }+ T( \hat{r}^s ) - T( r^s) \nonumber \\ 
&\stackrel{ \eqref{eq:bound_stage_cost_new}, \eqref{eq:bound_short_cost}, \eqref{eq:evo_offset}}{\leq }  
-c\|r^s-r^d_r(\mathbb{A})\|^2,
\end{align}
where $c>0$ is defined as follows
\begin{align*}
c:=&
 \left(\frac{1}{a_1} + \epsilon\gamma^2 k_1^{\ell}\right)\bar{a}  +
 \left(\epsilon\gamma a_2k_0^T \sqrt{\frac{\bar{a}}{a_1}}- k_1^T \right) \tau \\& + (\epsilon\gamma k_2^{\ell}k_0^{T,2}+\epsilon\gamma a_2 k_0^{T,2}) \tau^2   .
\end{align*}
The inequality above means that the candidate cost is smaller than the optimal one (i.e., satisfaction of \eqref{eq:desired_result_cost}), which results in a contradiction.
\hfill $\square$
\color{black}
%
%

\subsection{Recursive Feasibility of the MPC problem \eqref{eq:proposed_MPC}}
\label{sec:recursive_feasibility}
Given the optimal backup trajectory $(x^{B,*}_{\cdot| t}, u^{B,*}_{\cdot | t}),~ r^{sB,*}_{t}$ resulting from \eqref{eq:proposed_MPC} at time $t$, the candidate backup trajectory $(\hat{x}^{B}_{\cdot| t+1}, \hat{u}^{B}_{\cdot | t+1}), \hat{r}^{sB}_{t+1}$, and the candidate learning trajectory $(\hat{x}^{L}_{\cdot| t+1}, \hat{u}^{L}_{\cdot | t+1}), \hat{r}^{sL}_{t+1}$ feasible at time $t+1$ are both defined as follows
\begin{align}
\label{eq:candidate_backup}
\hat{r}^{sB}_{t+1} =\hat{r}^{sL}_{t+1}& := r^{sB,*}_{t}, \nonumber \\
\hat{x}^{B}_{0 |t+1}=\hat{x}^{L}_{0 |t+1} & := x_{t+1}, \nonumber \\
\hat{x}^{B}_{k-1|t+1}=\hat{x}^{L}_{k-1|t+1} & := x^{B,*}_{k| t}, \quad k=1, \dots, N, \nonumber\\
\hat{u}^{B}_{k-1|t+1}=\hat{u}^{L}_{k-1|t+1} & := u^{B,*}_{k| t}, \quad k=1, \dots, N-1, \nonumber\\
\hat{u}^{B}_{N-1| t}=\hat{u}^{L}_{N-1| t} &:= \kappa_f(\hat{x}^{B}_{N - 1|t+1}, \hat{r}^{sB}_{t+1}), \nonumber\\
\hat{x}^{B}_{N|t+1}=\hat{x}^{L}_{N|t+1} & := f(\hat{x}^{B}_{N-1|t+1}, \hat{u}^{B}_{N-1| t}).
\end{align}
Note that due to the potentially changing estimated constraint set $\mathbb{Z}_t$, it is not possible to ensure that the optimal learning trajectory $(x^{L,*}_{\cdot| t}, u^{L,*}_{\cdot | t}), r^{sL,*}_{t}$ computed at time $t$ can be employed for constructing the candidate learning trajectory at time $t+1$. For this reason, both the learning and backup trajectories are defined based on the optimal backup trajectory.

In the following, we analyse how the candidate trajectories in \eqref{eq:candidate_backup} satisfy the constraints \eqref{eq:MPC_backup_constraints}, \eqref{eq:MPC_terminal_backup}, \eqref{eq:MPC_learning_constraints}, \eqref{eq:MPC_terminal_learning}, and \eqref{eq:sum_risk_negative}. Note that the remaining constraints in \eqref{eq:proposed_MPC} are trivially satisfied, and therefore will not be discussed.

\paragraph*{Satisfaction of \eqref{eq:MPC_backup_constraints}, \eqref{eq:MPC_terminal_backup}, \eqref{eq:MPC_learning_constraints}, and \eqref{eq:MPC_terminal_learning}}
Satisfaction of these constraints is ensured based on the properties of the terminal region $\mathbb{X}_f$ shown in \eqref{eq:terminal_ingredients}, on the candidate trajectories defined in \eqref{eq:candidate_backup}, and on how the safe set $\mathbb{E}_t$ is updated, cf.
 \eqref{eq:set_update}.
\paragraph*{Satisfaction of \eqref{eq:sum_risk_negative}}
Firstly, note that based on satisfaction of \eqref{eq:sum_risk_negative} up to time $t$, and since we start with $S_0\geq0$, we have that the following holds
\begin{align}
\label{eq:lower_bound_S}
\underbrace{ S_t - F^{*}_{t}  + \hat{F}_t}_{\stackrel{\eqref{eq:function_S_t}}{=} S_{t+1}}\stackrel{\eqref{eq:sum_risk_negative}}{\geq} 0 \Leftrightarrow S_{t+1} \geq 0, \quad \forall t\geq 0.
\end{align}
Considering the candidate backup trajectory defined in \eqref{eq:candidate_backup}, we have that
\begin{align}
\label{eq:candidate_vs_previus}
& \epsilon V_N(\hat{x}_{\cdot|t+1}^{B}, ~\hat{u}_{\cdot|t+1}^{B},~ \hat{r}^{sB}_{t+1}) + T(\hat{r}_{t+1}^{sB}) \nonumber \\
\stackrel{\eqref{eq:terminal_ingredients}, \eqref{eq:candidate_backup} }{\leq} & \epsilon V_N(x_{\cdot|t}^{B,*}, ~u_{\cdot|t}^{B,*},~ r^{sB,*}_{t}) + T(r_{t}^{sB,*}) - \epsilon\ell(x^{B,*}_{t},u^{B,*}_{t}, r^{sB,*}_{t})\nonumber\\
\stackrel{\eqref{eq:F_optimal}}{=}& F^*_t - \epsilon\ell(x^{B,*}_{t},u^{B,*}_{t}, r^{sB,*}_{t})
\stackrel{\eqref{eq:F_backup}}{\leq} \hat{F}_{t+1} \leq \underbrace{S_{t+1}}_{\stackrel{\eqref{eq:lower_bound_S}}{\geq} 0} +\hat{F}_{t+1} ,
\end{align}
which ensures satisfaction of \eqref{eq:sum_risk_negative} at time $t+1$ as well.

\subsection{Convergence (Satisfaction of  \eqref{eq:manifold_theorem} and \eqref{eq:decrease_theorem})}
\label{sec:convergence}
For the sake of simplicity, instead of defining $S_t$ recursively as done in \eqref{eq:function_S_t}, we reformulate $S_t$ as follows
\begin{align}
\label{eq:reformulation_S_t}
S_{t+1} \stackrel{\eqref{eq:function_S_t}}{=} & S_t +\hat{F}_t -  F^*_{t} =  S_0 + \sum_{i=0}^t  \hat{F}_i-F^*_{i} \nonumber \\ 
= &  S_0 +  \hat{F}_0 - \hat{F}_{t+1} + \sum_{i=0}^t  \hat{F}_{i+1}-F^*_{i} \nonumber \\
\stackrel{\eqref{eq:F_backup}}{=} & S_0 +  \hat{F}_0 - \hat{F}_{t+1} - \sum_{i=0}^{t} \epsilon\alpha\ell(x^{B,*}_{i},u^{B,*}_{i}, r^{sB,*}_{i}).
\end{align}

Therefore, based on \eqref{eq:sum_risk_negative} and on \eqref{eq:reformulation_S_t}, the following holds
\begin{align}
\label{eq:bound_F_epsilon}
 F_{ t+1}^* \stackrel{\eqref{eq:sum_risk_negative}, \eqref{eq:F_optimal}}{\leq} & S_{t+1} + \hat{F}_{t+1} \nonumber \\
\stackrel{\eqref{eq:reformulation_S_t}}{=} & S_0 - \hat{F}_{t+1} + \hat{F}_{0} - \sum_{i=0}^{t} \epsilon\alpha\ell(x^{B,*}_{i},u^{B,*}_{i}, r^{sB,*}_{i}) +\hat{F}_{t+1} \nonumber \\
= & S_0 + \hat{F}_{0} - \sum_{i=0}^t \epsilon\alpha\ell(x^{B,*}_{i},u^{B,*}_{i}, r^{sB,*}_{i}).
\end{align}
By re-arranging \eqref{eq:bound_F_epsilon}, and taking the infinite average sum of both sides, we have
\begin{align}
\label{eq:convergence_to_steady_state}
& \limsup_{T\rightarrow \infty} \frac{F_{ T+1}^* - \hat{F}_{0} + \sum_{t=0}^T  \epsilon\alpha \ell(x^{B,*}_{i},u^{B,*}_{i}, r^{sB,*}_{i})}{T}  \nonumber \\
\stackrel{\eqref{eq:bound_F_epsilon}}{\leq }& \limsup_{T\rightarrow \infty} \frac{ S_0 }{T} = 0 \nonumber \\
\Rightarrow & \lim_{T\rightarrow \infty} \frac{ \sum_{t=0}^T  \epsilon\alpha \ell(x^{B,*}_{t},u^{B,*}_{t}, r^{sB,*}_{t})}{T} \leq 0 \nonumber \\
\Rightarrow & \lim_{T\rightarrow \infty} \ell(x^{B,*}_{t},u^{B,*}_{t}, r^{sB,*}_{t}) = 0 \nonumber \\
\stackrel{\text{Ass. \ref{ass:bound_cost_function}}}{\Rightarrow}& \lim_{T\rightarrow \infty} V_N(x_{\cdot|t}^{B,*}, ~u_{\cdot|t}^{B,*},~ r^{sB,*}_t) = 0.
\end{align}
The result above holds since $F_{ t}^*$ is lower and upper bounded for all $t\geq0$. In particular, the lower bound for $F_{t}^*$ (i.e., $F_{t}^*\geq 0$) can be trivially ensured by construction of the cost function in \eqref{eq:cost_function}, and thanks to the lower and upper bounds of the stage costs discussed in Assumption \ref{ass:stage_cost}, while the upper-bound is shown in \eqref{eq:bound_F_epsilon}, based on the fact that both $S_0$ and $\hat{F}_{0}$ are bounded.
Based on both \eqref{eq:sum_risk_negative} and on \eqref{eq:convergence_to_steady_state}, then it is trivial to show satisfaction of \eqref{eq:decrease_theorem}.

Note that in \eqref{eq:convergence_to_steady_state} we only showed that the system converges to a feasible artificial setpoint $r^{sB,*}_{\infty} \in \mathbb{S}(\mathbb{Z})$, i.e.,
\begin{align}
\label{eq:proof_convergence}
\lim_{t\rightarrow \infty} z_t = r^{sB,*}_{\infty} \in \mathbb{S}(\mathbb{Z}), \quad z_t=(x_t, u_t).
\end{align}
However, in order to show satisfaction of \eqref{eq:manifold_theorem}, we additionally need to show that $r^{sB,*}_{\infty} \notin\mathbb{M}(\epsilon,N,\mathbb{H}(r^{s,*}_{\infty}))$, and $r^{sB,*}_{\infty} \in \text{bound}(\mathbb{S}(\mathbb{Z}))$.

\paragraph*{$r^{sB,*}_{\infty} \notin \mathbb{M}(\epsilon,N,\mathbb{H}(r^{s,*}_{\infty}))$}  The following part will be done by contradiction and therefore we assume that the system converges to a transitory setpoint $r^{sB,*}_{\infty} \in \mathbb{M}(\epsilon,N,\mathbb{H}(r^{s,*}_{\infty}))$.

Firstly, based on \eqref{eq:convergence_to_steady_state}, we know that there exists a time $t_{\delta}\geq0$ such that 
\begin{align}
\label{eq:small_stage_cost}
V_N(x_{\cdot|t_{\delta}}^{B,*}, u_{\cdot|t_{\delta}}^{B,*},r^{sB,*}_{\infty})  \leq \delta,
\end{align}
for any $\delta>0$.
According to Definition \ref{def:transitory_steady_state}, we know that since $ r^{sB,*}_{\infty}$ is a transitory artificial setpoint, then there exists a different artificial setpoint $\bar{r}^{sB}$ with $T(\bar{r}^{sB})< T(r^{sB,*}_{\infty})$ such that the following holds
\begin{align}
\label{eq:better_steady_state}
& \epsilon\left[\ell(x_{0|t_{\delta}}^{B,*},u_{0|t_{\delta}}^{B,*}, \bar{r}^{sB}) + V^*_{N-1}(x_{1|t_{\delta}}^{B,*}, \bar{r}^{sB}, \mathbb{M}(\epsilon,N,\mathbb{H}(r^{s,*}_{\infty})))\right]  \nonumber \\
& +T(\bar{r}^{sB})<  \epsilon V_N^*(x_{t_{\delta}},r^{sB,*}_{\infty},\mathbb{M}(\epsilon,N,\mathbb{H}(r^{s,*}_{\infty}))) + T(r^{sB,*}_{\infty}).
\end{align}
The new candidate solution related to the setpoint $\bar{r}^{sB}$ is feasible for the MPC scheme \eqref{eq:proof_convergence} since it    trivially satisfies \eqref{eq:sum_risk_negative} based on \eqref{eq:better_steady_state} and \eqref{eq:MPC_first_input}.
This implies that the overall candidate cost function satisfies the following
\begin{align}
\label{eq:optimal_cost_contradiction}
& V_N(x^{L*}_{\cdot|t},u^{L*}_{\cdot|t},r^{sL,*}_{t}) +T(r^{sL,*}_t)+ \epsilon T(\bar{r}^{sB}) \nonumber  \\
< & V_N(x^{L*}_{\cdot|t},u^{L*}_{\cdot|t},r^{sL,*}_{t}) +T(r^{sL,*}_t)+ \epsilon T(r^{sB,*}_{\infty}).
\end{align}
Hence, we showed that there exists a feasible candidate trajectory that has a smaller overall cost function, which, therefore, results in a contradiction.

\paragraph*{$r^{sB,*}_{\infty} \in \text{bound}(\mathbb{S}(\mathbb{Z}))$} This part is again done by contradiction, and by making use of Proposition \ref{prop:sliding_steady_state}, as well as the definition of the set $\mathbb{H}$ in Assumption \ref{ass:neighborhood}. In the following, we assume that $r^{sB,*}_{\infty}\in \text{int}(\mathbb{S}(\mathbb{Z}))$. Based on \eqref{eq:proof_convergence}, we know that 
there exists a time $t^*$ such that $z_t=(x_t, u_t)$ satisfies
\begin{align*}
\|z_t - r^{sB,*}_{\infty} \|^2 < h-\lambda, \quad \forall t\geq t^*,
\end{align*}
which, based on the construction of the set $\mathbb{H}$ in Assumption \ref{ass:neighborhood} and on the definition of the manifold in \eqref{eq:manifold}, it implies that $r^{sB,*}_{\infty} \in \text{int}(\mathbb{H}(z_t))$, $r^{sB,*}_{\infty} \in \text{int}(\mathbb{S}(\mathbb{H}(z_t)))$, for all $t^*\geq 0$. Since we already showed that $r^{sB,*}_{\infty}$ is not a transitory setpoint, then based on Proposition \ref{prop:sliding_steady_state}, we also know that $r^{sB,*}_{\infty} \in \text{bound}(\mathbb{S}(\mathbb{H}(z_t)))$. Due to the contradiction obtained above, one can show that the only possible solution is that $r^{sB,*}_{\infty}\in \text{bound}(\mathbb{S}(\mathbb{Z}))$.
%
%
%

%% file: example.tex
\begin{figure*}[t]
   \centering
\begin{tabular}{cc}
\includegraphics[width=0.45\textwidth]{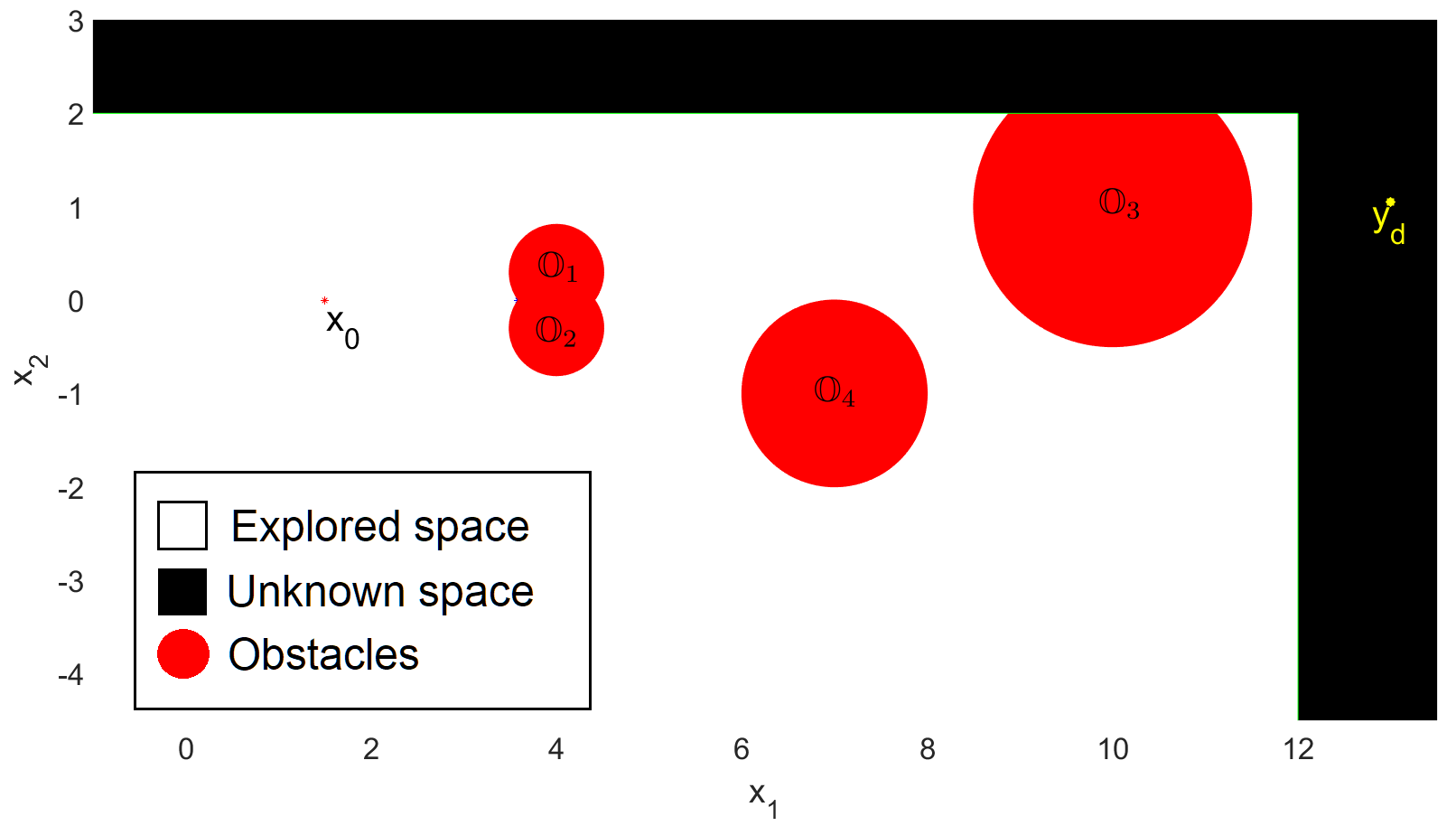}
&
\includegraphics[width=0.45\textwidth]{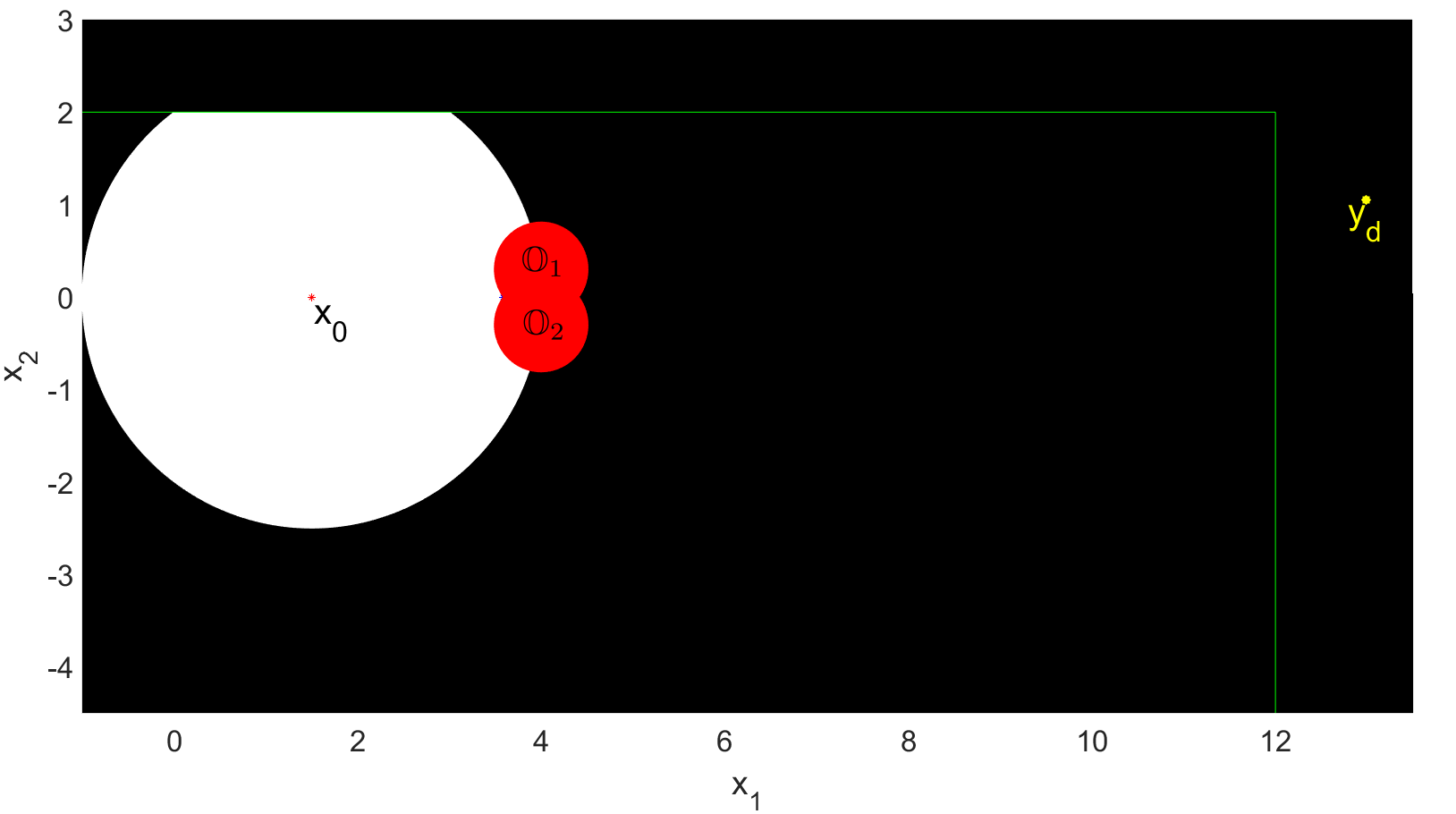}\\
(a) Ideal configuration where the entire set $\mathbb{Z}$ & (b) Actual configuration at $t=0$\\(including obstacles) is known, i.e., $\mathbb{E}_0=\mathbb{Z}$ &  where only a neighborhood of $x_0$ is known\\
\hspace{-0.4cm}\includegraphics[width=0.45\textwidth]{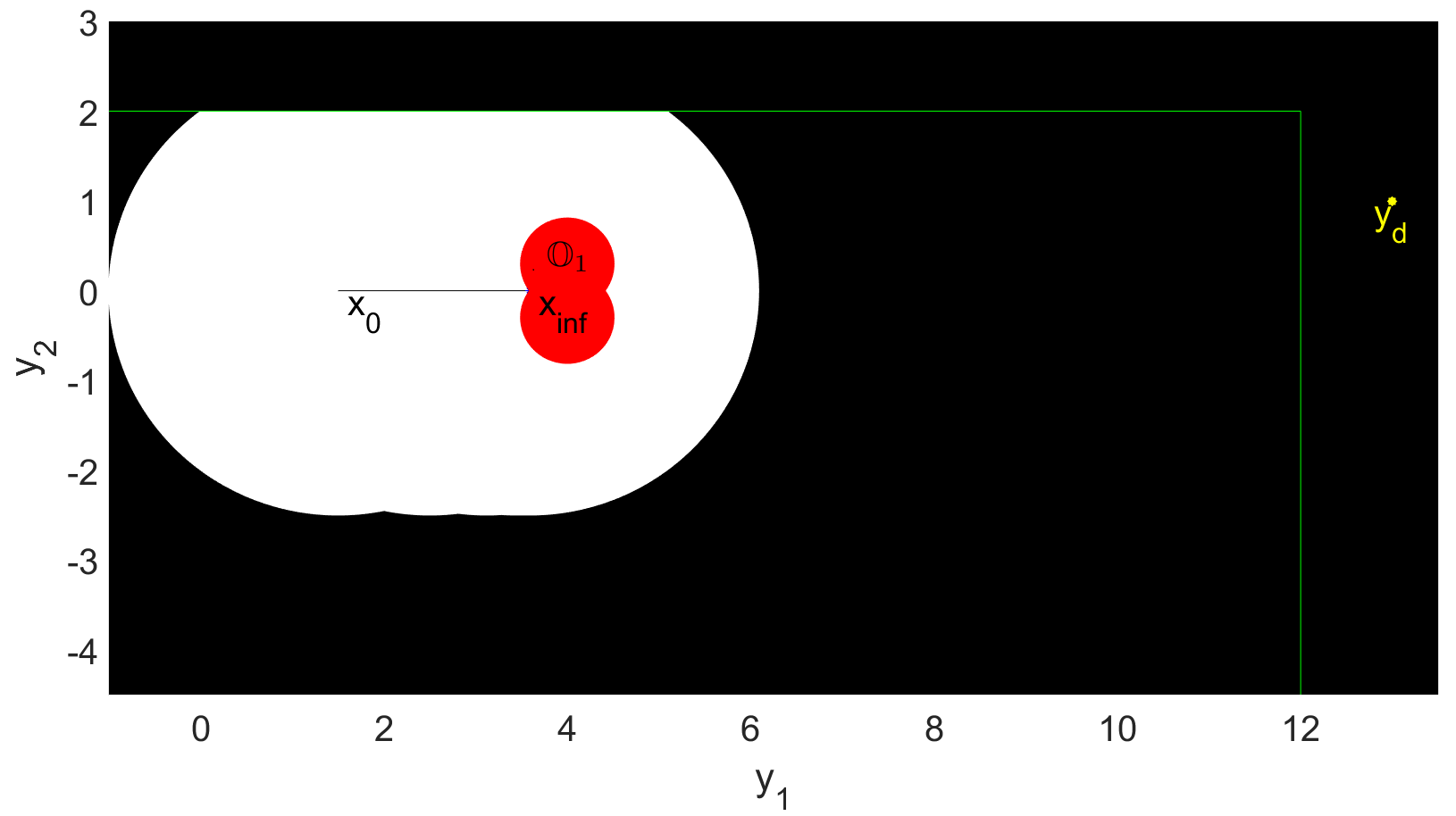}&
\includegraphics[width=0.45\textwidth]{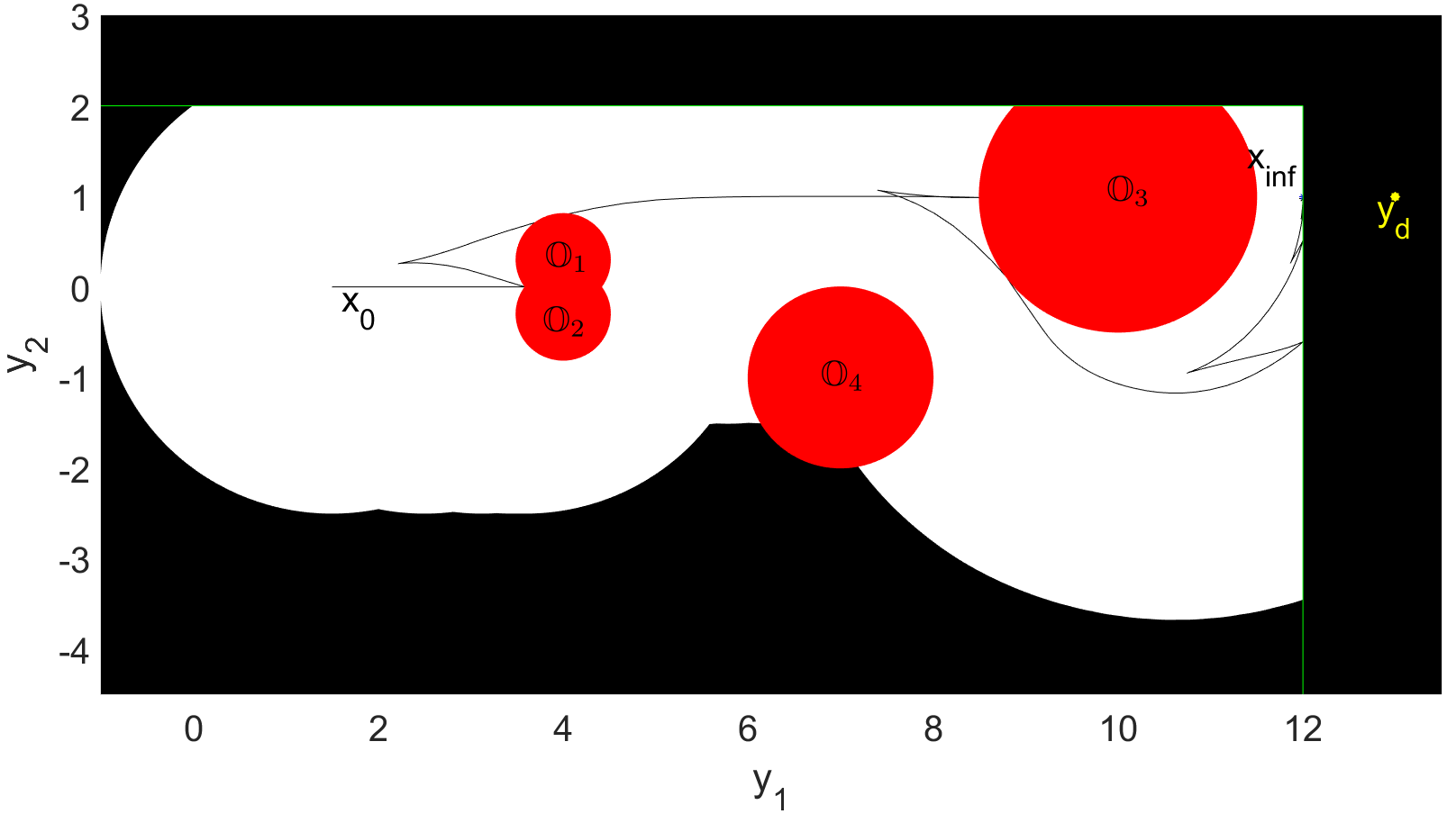}\\
(c) Closed-loop of system \eqref{eq:system_example} controlled & (d) Closed-loop of system \eqref{eq:system_example} controlled \\
 with \eqref{eq:proposed_MPC} without \eqref{eq:sum_risk_negative} & with \eqref{eq:proposed_MPC} with \eqref{eq:sum_risk_negative}
\end{tabular}
   \caption{\RS{Evolution of system \eqref{eq:system_example} with and without constraint \eqref{eq:sum_risk_negative}.}}
\label{fig:example}
\end{figure*}

\section{Example}
\label{sec:example}
In this example, we consider a car traveling in a partially unknown space where several unknown obstacles are present, with the goal to get as close as possible to the unreachable output $y^d$, as shown in Figure \ref{fig:example}. 
We compare the proposed approach with the results shown in \cite{limon2018nonlinear} and extended in Section \ref{sec:preliminaries}, as well as with the MPC scheme \eqref{eq:proposed_MPC} without the implementation of condition \eqref{eq:sum_risk_negative}. We want to additionally emphasize that the computationally tractable approach for constructing convex sets proposed in \cite{tordesillas2021faster} can be combined with the proposed MPC framework, in order to reduce its computational complexity, if necessary.

\paragraph*{Problem setup}
We consider a nonlinear system representing an autonomous car, defined as follows
\begin{align}
\label{eq:system_example}
\dot{x}=\begin{bmatrix}
\dot{x_1}\\
\dot{x_2}\\
\dot{\psi}\\
\dot{v}\\
\dot{\beta}
\end{bmatrix} = 
\begin{bmatrix}
v\cdot \cos(\psi+\beta)\\
   v\cdot \sin(\psi+\beta)\\
   \frac{v}{l_r}\cdot \sin(\beta)\\
   u_1\\
   u_2
\end{bmatrix}, \quad 
y = \begin{bmatrix}
y_1\\ y_2
\end{bmatrix}=\begin{bmatrix}
x_1\\ x_2
\end{bmatrix},
\end{align}
where $x_1\in [-\infty, 12]$ and $x_2\in [-\infty, 2]$ represent the coordinates of the center of mass of the car, $\psi\in [-\infty, \infty]$ indicates its rotation, $v\in [-1m/s, 36m/s]$ is the longitudinal velocity, and $\beta \in [-37^\circ, 37^\circ]$ is the angle of the current velocity of the center of mass with respect to the longitudinal axis of the car, compare \cite{soloperto2019collision}. The inputs $u_1\in [-10m/s^2, 1m/s^2]$, and $u_2\in [-10^\circ/s, 10^\circ/s]$ control the acceleration of the car and the velocity of the steering angle, respectively. We define the set $\mathbb{Z}_{sc}$ as the hyper-box that encodes the state and input constraints introduced above.
The parameter $l_r=1.7m$ represents the distance from the center of the mass of the car to the rear axle. 
For the implementation we use a Euler discretization with a sampling time of $h=200$ ms, as done in \cite{kong2015kinematic}.
The control goal is to get as close as possible to the desired unreachable output $y_d = [12, 1]$.

\paragraph*{Obstacles description}
We define the obstacles as ellipsoidal sets in the output space $y=(x_1,x_2)$, compare Figure \ref{fig:example}. Each obstacle is formally defined as $\mathbb{O}_i:=\{x\in \mathbb{R}^{n+m}| \|[x_1,x_2]^{\top}-o_i\|_{I_2}\leq r_i\}$, where $o_i\in \mathbb{R}^2$ is the center of the obstacle, $r_i\in \mathbb{R}^2$ its radius, and, $I_2$ the identity matrix, with $i=1, \dots, 4$.
The values of $o_i$ and $r_i$ are chosen as follows
\begin{align*}
o_1 = \begin{bmatrix}
4.0\\ 0.3
\end{bmatrix},~o_2 = \begin{bmatrix}
4.0\\ -0.3
\end{bmatrix},~o_3 = \begin{bmatrix}
10.0\\ 1.0
\end{bmatrix},~o_4 = \begin{bmatrix}
7.0\\ -1.0
\end{bmatrix},\\
r_1 = 0.51,~r_2 = 0.51,~r_3 = 1.50,~r_4 = 1.00.
\end{align*}
Note that each  obstacle is  visible to the vehicle if and only if it is in the \textit{neighborhood} $\mathbb{H}(z_t)$ of the current output $y_t$ of the system, where $\mathbb{H}$ is introduced below.

The constraint set $\mathbb{Z}$ is defined based on the set $\mathbb{Z}_{sc}$, and the obstacles $\mathbb{O}_i$ described above, i.e.,
\begin{align*}
\mathbb{Z}:= \mathbb{Z}_{sc} \cap(\cup_{i=1}^4 \mathbb{O}_i).
\end{align*}

\paragraph*{Set $\mathbb{H}(z_t)$ and update of $\mathbb{E}_t$ and $\mathbb{Z}_t$}
In this example, we consider the case where the vehicle is equipped with a sensor (e.g., a LIDAR sensor) that is able to perceive the surrounding $\mathbb{H}(z_t)$ of the vehicle, up to a radius of $2.5$m in the output space
\begin{align*}
\mathbb{H}(z_t) := \{x\in \mathbb{R}^{n+m}| \|[x_1,x_2]^{\top}-y_t\|_{I_2}\leq 2.5\}\cap \mathbb{Z}.
\end{align*}
Based on the set $\mathbb{H}(z_t)$, we define, at each time $t\geq0$, the set $\mathbb{E}_t$ as follows
\begin{align*}
\mathbb{E}_t := \mathbb{Z}_{sc}\cap( \cup_{i=0}^t \mathbb{H}(z_i)) = \mathbb{E}_{t-1} \cup \mathbb{H}(z_t).
\end{align*}
On the contrary, the set $\mathbb{Z}_t$ is initialized as $\mathbb{Z}_0:= \mathbb{Z}_{sc}$, and additionally includes the obstacles $\mathbb{O}_i$ only if they are in a neighborhood $\mathbb{H}(z_t)$ of the system itself. This implies that the learning trajectory can freely explore the unknown space $\mathbb{Z}_t\backslash\mathbb{E}_t$, while the backup trajectory is forced to stay inside the explored set $\mathbb{E}_t$.
\paragraph*{MPC scheme}
We use a standard quadratic stage cost with weighting matrices $Q:= \text{diag}[1, 1, 10^{-5}, 10^{-5}, 10^{-5}]$ and $R:=\text{diag}[1, 1]$. Additionally, to simplify the implementation, we consider a terminal equality constraints w.r.t. the artificial setpoints $r^s$. 
The offset cost is defined as $T(y) = \|y-y^d\|^2$, while we choose $\epsilon=0.01$, $\alpha=1$ and an MPC horizon of $N=50$, which corresponds to $100s$ in the given setup. Finally, the initial state is set as 
\begin{align*}
x_0:=[1.5, 0, 0, 5, 0]^{\top},
\end{align*}
as shown in Figure \ref{fig:example}.
\paragraph*{Discussion}
In Figure \ref{fig:example}-(a), we show the ideal configuration where the entire set $\mathbb{Z}$ is known, i.e., all the obstacles are visible at time $t=0$. On the contrary, in Figure \ref{fig:example}-(b), we  show the actual configuration that we employ in the considered example at $t=0$. In particular, we see that only the neighborhood of the initial set $x_0$ is visible, which includes also the obstacles $\mathbb{O}_1$ and $\mathbb{O}_2$, together with the constraint $x_2\leq 2$. In Figure \ref{fig:example}-(c), we show how the MPC scheme \eqref{eq:proposed_MPC} without the constraint \eqref{eq:sum_risk_negative} behaves. Specifically, we see that the system starts in $x_0$ and converges to the artificial steady state $[x_1,x_2]^{\top}\approx[3.5, 0]^{\top}$, which corresponds to a non-transitory setpoint due to the presence of obstacles. This happens because the vehicle is not able to get closer to the desired output $y_d$ since this would require a phase where it firstly moves farther away from the desired output $y^d$, and only after a proper maneuver it is able to overcome the obstacles and it can then get closer to $y^d$. Note that Figure \ref{fig:example}-(c) also represents the behavior of system \eqref{eq:system_example} under the MPC scheme \eqref{eq:standard_MPC_scheme}, which is analogous to the approach proposed in \cite{limon2018nonlinear}.

Finally, in Figure \ref{fig:example}-(d), we show the closed-loop system of system \eqref{eq:system_example} under the control law resulting from \eqref{eq:proposed_MPC} when \eqref{eq:sum_risk_negative} is included. In this case, we see that even though the system reaches several non-transitory setpoints, it is nevertheless able to move backward (i.e., temporarily increasing the overall cost) and then move towards a different artificial setpoint that is closer to the desired output $y^d$. At time $t=60s$, the system converges to the output $y = [12, 1]$, which corresponds to the closest reachable output $y^d_r(\mathbb{Z})$, as defined in \eqref{eq:best_reachable_setpoint}.

In conclusion, this example illustrates the capabilities of the proposed MPC scheme to escape from local minima, while ensuring safety (i.e., closed-loop constraint satisfaction),  convergence, and the possibility to consider a non-reachable output $y^d$.